\documentclass[11pt,leqno]{article}

\usepackage{amsmath,amsfonts,amscd,amssymb,theorem}

\long\def\comment#1\endcomment{}

\comment
\endcomment

\makeatletter
\begingroup
\gdef\th@dotted{\normalfont\itshape
  \def\@begintheorem##1##2{%
        \item[\hskip\labelsep \theorem@headerfont ##1\ ##2.]}%
\def\@opargbegintheorem##1##2##3{%
   \item[\hskip\labelsep \theorem@headerfont ##1\ ##2\ (##3).]}}
\endgroup
\makeatother

\theoremstyle{dotted}

\newtheorem{theorem}{Theorem}[section]
\newtheorem{lemma}[theorem]{Lemma}
\newtheorem{conj}[theorem]{Conjecture}
\newtheorem{prop}[theorem]{Proposition}
\newtheorem{corr}[theorem]{Corollary}

\makeatletter
\begingroup
\gdef\th@upshape{\normalfont
  \def\@begintheorem##1##2{%
        \item[\hskip\labelsep \theorem@headerfont ##1\ ##2.]}%
\def\@opargbegintheorem##1##2##3{%
   \item[\hskip\labelsep \theorem@headerfont ##1\ ##2\ (##3).]}}
\endgroup
\makeatother

\theoremstyle{upshape}

\newtheorem{defn}[theorem]{Definition}
\newtheorem{remark}[theorem]{Remark}
\newtheorem{example}{Example}[section]

\makeatletter
\renewcommand{\subsection}{\@startsection{subsection}{2}{0pt}{-3ex
plus -1ex minus -0.2ex}{-2mm plus -0pt minus
-2pt}{\normalfont\bfseries}} \makeatother

\makeatletter
\@addtoreset{equation}{section}
\makeatother

\newcommand{\cntrct}                
{\hspace{2pt}\raisebox{1pt}{\text{$\lrcorner$}}\hspace{2pt}}

\newcommand{\proof}[1][Proof.]{\smallskip\noindent{\em #1}}
\def\endproof{\hfill\ensuremath{\square}\par\medskip}

\newcounter{stepp}
\newcommand{\step}{\medskip\stepcounter{stepp}

\noindent
{\em Step \thestepp.\/}}

\def\eqref#1{\thetag{\ref{#1}}}

\let\latexref=\ref
\def\ref#1{{\normalfont{\latexref{#1}}}}

\newcommand{\wt}{\widetilde}
\newcommand{\wh}{\widehat}
\newcommand{\whtimes}{\widehat{\times}}

\newcommand{\ratto}{\dasharrow}        

\newcommand{\idot}{{\:\raisebox{1pt}{\text{\circle*{1.5}}}}}
\newcommand{\hdot}{{\:\raisebox{3pt}{\text{\circle*{1.5}}}}}
\newcommand{\g}{{\mathfrak{g}}}
\newcommand{\h}{{\mathfrak{h}}}

\newcommand{\tw}{ {(1)} }
\newcommand{\Fr}{{\sf Fr}}
\newcommand{\ad}{\operatorname{\sf ad}}
\newcommand{\cchar}{\operatorname{\sf char}}
\newcommand{\Pic}{\operatorname{\sf Pic}}
\newcommand{\Br}{\operatorname{\sf Br}}

\newcommand{\Spec}{\operatorname{Spec}}

\newcommand{\sspec}{\operatorname{{\cal S}{\it pec}}}

\newcommand{\Shv}{\operatorname{Shv}}
\newcommand{\fmod}{{\amod^{\text{{\tt\tiny fg}}}}}
\newcommand{\amod}{{\text{\rm -mod}}}

\newcommand{\C}{{\mathbb C}}
\newcommand{\R}{{\mathbb R}}

\newcommand{\Z}{{\mathbb Z}}

\newcommand{\A}{{\mathbb A}}

\newcommand{\cala}{{\mathcal A}}
\newcommand{\calr}{{\mathcal R}}

\newcommand{\barA}{\overline{A}}
\newcommand{\barD}{\overline{D}}

\newcommand{\Hh}{{\mathbb H}}
\newcommand{\HH}{{\mathcal H}}

\newcommand{\Pp}{{\mathbb P}}
\newcommand{\PP}{{\bf P}}

\newcommand{\X}{{\mathfrak X}}

\newcommand{\D}{{\cal D}}

\newcommand{\F}{{\cal F}}
\newcommand{\M}{{\cal M}}

\newcommand{\LL}{{\cal L}}
\newcommand{\T}{{\cal T}}

\newcommand{\Y}{{\cal Y}}

\newcommand{\E}{{\cal E}}
\newcommand{\K}{{\cal K}}

\newcommand{\G}{{\mathcal G}}

\newcommand{\calo}{{\cal O}}

\newcommand{\m}{{\mathfrak m}}

\newcommand{\gm}{{\mathbb{G}_m}}

\renewcommand{\phi}{\varphi}

\renewcommand{\dim}{\operatorname{\sf dim}}
\newcommand{\codim}{\operatorname{\sf codim}}

\newcommand{\gr}{\operatorname{\sf gr}}

\newcommand{\Id}{\operatorname{{\sf Id}}}

\newcommand{\Ext}{\operatorname{Ext}}

\newcommand{\End}{\operatorname{End}}
\newcommand{\eend}{\operatorname{{\cal E}{\it nd}}}
\newcommand{\Aut}{\operatorname{{\sf Aut}}}
\newcommand{\Der}{\operatorname{{\sf Der}}}
\newcommand{\Loc}{\operatorname{{\sf Loc}}}

\newcommand{\Rhom}{\operatorname{{\bf R}{\cal H}{\it om}}}

\newcommand{\Maps}{\operatorname{{\sf Maps}}}

\newcommand{\Supp}{\operatorname{\sf Supp}}

\title{Geometry and topology of symplectic resolutions}

\author{D. Kaledin\thanks{Partially supported by CRDF grant
RUM1-2694.}}

\begin{document}

\maketitle

\tableofcontents

\section*{Introduction.}

By Hironaka, every singular algebraic variety $Y$ over $\C$ admits a
resolution of singularities -- that is, a smooth algebraic variety
$X$ equipped with a projective birational map $X \to Y$. In many
problems of algebraic geometry the mere existence of $X$ is enough,
but sometimes it is not. Especially when algebraic geometry is being
used as a tool in some other area of mathematics, more control over
the resolution $X$ is needed. 

This is very prominently the case, for instance, in Geometric
Representation Theory (see e.g. \cite{GC}). Ideally, given a
singular variety $Y$ encoding some representation theory problem,
one wants a resolution $X$ which is semismall (that is, $\dim X
\times_Y X = \dim X$), and with some restrictions on the topology of
the fibers. If $Y$ carries some group action, one also wants $X$ to
be equivariant with respect to this action.

In many cases it is actually possible to achieve the ideal. For
instance, if $Y$ is the nilpotent cone in the adjoint representation
of a semisimple Lie group, then it admits a semismall resolution $X$
known as {\em Springer resolution}. The resolution is equivariant
with respect to all possible group actions on $Y$. Its fibers,
although singular, cohomologically behave in the same way as smooth
homogeneous spaces: all the cohomology groups are pure with respect
to the weight filtration, and are in fact spanned by classes of
algebraic cycles. A completely parallel picture holds for the
so-called {\em quiver varieties} of H. Nakajima, and for Hilbert
schemes of $n$ points on $\C^2$. Some additional pieces of structure
are present in all these cases on the resolution $X$, in particular,
$X$ is equipped with a holomorphic symplectic from.

The standard proofs of these facts (e.g. \cite{DLP}) work by
explicit constructions and rely heavily on the specific geometry of
the variety $Y$ in question.

\bigskip

The goal of the present paper is to report on a recent series of
results that somewhat changes this conventional perspective. Namely,
upon closer inspection, it turns out that the holomorphic symplectic
form, an auxiliary and almost accidental piece of structure on the
resolution $X$, actually insures all of the other good properties it
enjoys -- the semismall property, the cohomological purity of the
fibers, and so on and so forth. Moreover, the theory can be pushed
through so far as to give a complete algebraic description of the
derived category of coherent sheaves on the $X$. This gives new
information even in the well-studied cases such as the Springer
resolution or the Hilbert scheme.

\bigskip

Since the only thing needed from $X$ is the holomorphic symplectic
form, the results we are going to describe lie entirely within
Algebraic Geometry (or even Algebraic Symplectic Geometry, if such a
thing exists at present). Thus no knowledge of Geometric
Representation Theory is needed, nor assumed. Moreover, while most
known applications come from geometric representation theory, all
the results can also be used in the local study of contractions of
compact holomorphic symplectic and hyperk\"ahler manifolds -- or,
more generally, in that part of the Minimal Model Program which
deals with varieties with trivial canonical bundle. In particular,
some of the results on the derived category are actually difficult
conjectures which should hold in larger generality, and at least for
general Calabi-Yau varieties (see e.g. the general program sketched
in \cite{BO1}, \cite{BO2}). In the holomorphic symplectic case,
these conjectures can actually be proved.

We should warn the reader that the scope of this paper is limited --
we essentially restrict ourselves to giving an overview of the
papers \cite{ka1}, \cite{bk1}, \cite{bk2}, \cite{bk3} and
\cite{k1}. We do {\em not} attempt to give a general overview of
symplectic singularities, and we do not even mention a lot of fine
work -- a most notable omission is a series of papers \cite{N1},
\cite{N2}, \cite{N3}, \cite{N4} by Y. Namikawa and a paper \cite{FN}
by B. Fu and Y. Namikawa. In addition, we aim to be understandable
and brief, even at the cost of being precise. Thus some of the
proofs are omitted, and the other ones are only sketched. We always
give a precise references to original papers, which the reader who
wishes to see a complete proof should definitely consult. In
remarks, we allow ourselves even more imprecision, and the entire
last Subsection~\ref{quantu} should be treated as an extended
remark.

\subsection*{Acknowledgments.} The results reviewed in this paper
have been obtained through a long research project; the original
motivation for this project came from R. Bezrukavnikov, and part of
the research is joint work. It goes without saying that his help was
invaluable even in those parts which do not directly bear his
name. I am very grateful to the organizers of the Seattle meeting
for inviting me, and for generously allowing me several talk slots
instead of one, thus giving me an opportunity to present the results
in great detail.

\section{Definitions and general results.}\label{geo.sec}

Fix a base field $K$ of characteristic $0$. A convenient starting
point is the following definition introduced by A. Beauville
\cite{beau}.

\begin{defn}\label{s.sing.def}
A \!{\em symplectic singularity} is a normal irreducible algebraic
variety $Y$ over $K$ equipped with a non-degenerate symplectic form
$\Omega \in H^0(Y^{sm},\Omega^2)$ on the smooth locus $Y^{sm}
\subset Y$ which extends to a possibly degenerate symplectic form on
a smooth projective resolution $X \to Y$.
\end{defn}

Here and elsewhere in the paper {\em symplectic} is understood in
the algebraic sense -- it should not be confused with $C^\infty$
symplectic forms which appear in K\"ahler geometry. By a {\em
resolution} we understand a smooth variety $X$ equipped with a
projective birational map $X \to Y$. Originally, Beauville only
required the existence of the form $\Omega$; we prefer to include it
into the definition as a part of the data.

An easy observation (\cite{beau}) is that if $\Omega$ extends to one
smooth resolution $X$, it also extends to any other resolution $X'$
-- thus in Definition~\ref{s.sing.def} we can replace "a resolution"
with "any resolution" without any loss of generality.

In the present paper, we are mostly concerned with local study of
symplectic singularities -- in particular, $Y$ will be usually
assumed to be affine. Since we also assume $Y$ normal, we must have
$Y=\Spec H^0(X,\calo_X)$, so that, once a resolution $X$ is given,
it is no longer necessary to specify $Y$.

\medskip

Let us give some examples of symplectic singularities.

\begin{example}\label{dim2.exa}
$Y=W/G$, where $W$ is the $2$-dimensional vector space considered as
  an affine variety, and $G \subset SL(W)$ is a finite subgroup. In
  this classic case (see e.g. \cite{lau}), there exists a unique
  smooth resolution $X$ with trivial canonical bundle -- since we
  are in $\dim 2$, this is equivalent to having a symplectic form.
\end{example}

\begin{example}
$Y=\A^{2n}/S_n$, the quotient of the affine space of dimension $2n$
  by the symmetric group on $n$ letters -- equivalently, $Y$ is the
  $n$-th symmetric power of the affine plane $\A^2$. $X$ is the
  Hilbert scheme of $0$-dimensional subschemes of length $n$ in
  $\A^2$ (abbreviated to ``Hilbert scheme of $n$ points'').
\end{example}

\begin{example}\label{wreath.exa}
A combination of the previous two examples: $Y$ is the $n$-th
symmetric power of a quotient $Y_0=W/G$, $\dim W=2$, $G \subset
SL(W)$, $X$ is the Hilbert scheme of $n$ points on the canonical
symplectic resolution $X_0$ of $Y_0$.
\end{example}

\begin{example}\label{quo.exa}
$Y=V/G$ is the quotient of a symplectic vector space $V$ by a finite
  subgroup $G \subset Sp(V)$, $X$ is any resolution.
\end{example}

\begin{example}\label{spri.exa}
$Y \subset \g$ is the nilpotent cone in a Lie algebra
  $\g$ of a semisimple Lie group $G$, $X=T^*(G/B)$ is the cotangent
  bundle of the flag variety $G/B$ associated to $G$ (the {\em
  Springer resolution}).
\end{example}

\begin{example}\label{homo.exa}
More generally, $X=T^*(G/P)$ is the cotangent bundle to a
homogeneous variety $G/P$ associated to a parabolic subgroup $P
\subset G$ in a semisimple Lie group. $Y=\Spec H^0(X,\calo_X)$ is in
this case a closure of a certain nilpotent orbit in $\g$.
\end{example}

\begin{example}\label{nilp.exa}
Even more generally, $Y$ is the normalization of the closure of a
nilpotent orbit in a semisimple Lie algebra $\g$, $X$ is any
resolution.
\end{example}

\begin{example}\label{qui.exa}
$Y$ and $X$ are {\em quiver varieties} constructed by H. Nakajima
  \cite{naka} starting from certain combinatorial data.
\end{example}

We note that Examples~\ref{dim2.exa}-\ref{wreath.exa} are particular
cases of Example~\ref{quo.exa}, and
Examples~\ref{spri.exa}-\ref{homo.exa} are particular cases of
Example~\ref{nilp.exa}. The reason we have separated them from the
rest is that they actually satisfy a stronger assumption: there
exists a resolution $X$ to which the symplectic form $\Omega$
extends {\em as a non-degenerate $2$-form}. This property, unlike
the general definition of the symplectic singularity, depends on the
resolution. To emphasize this, we introduce the following
definition.

\begin{defn}\label{reso.defn}
A {\em symplectic resolution} is a smooth algebraic variety $X$ over
$K$ equipped with a closed non-degenerate $2$-form $\Omega$ such
that the canonical map $X \to Y=\Spec H^0(X,\calo_X)$ is a
birational projective map.
\end{defn}

Not all symplectic singularities have symplectic resolutions. In the
quotient singularity case (Example~\ref{quo.exa}), the only known
examples of symplectic resolutions are those in
Example~\ref{wreath.exa}. Moreover, it has been proved by
M. Verbitsky \cite{verb} that the existence of a symplectic
resolution yields a strong necessary condition on the subgroup $G
\subset Sp(V)$, and even this condition is not sufficient
(\cite{gk}). In the nilpotent orbit case (Example~\ref{nilp.exa}),
the existence question has been studied exhaustively by Baohua Fu
\cite{fu}; all the symplectic resolutions that do exist are covered
by Example~\ref{homo.exa}. Finally, in the quiver variety case a
resolution always exists, but this is in fact a corollary of a
certain technical assumption on the combinatorial data imposed in
\cite{naka}. If one drops this assumption, one obtains some quiver
varieties which do not admit a resolution (for instance, those
considered in \cite{kls}). Even more surprisingly, there are quiver
varieties which do admit a symplectic resolution, but {\em not} of
the quiver type -- such is the O'Grady singularity studied in
\cite{kl}.

An obvious source of smooth symplectic varieties is the cotangent
bundles, $X = T^*M$ for a smooth algebraic variety $M$. However,
this is not very promising from the point of view of symplectic
singularities. Indeed, there is the following conjecture variously
attributed to J.-P. Demailly, F. Campana, Th. Peternell, which is
very difficult, but kind of old and well-established.

\begin{conj}\label{dema}
Assume given a smooth algebraic variety $M$, let $X = T^*M$, and
assume that the natural map $X \to Y = H^0(X,\calo_X)$ is a
projective birational map. Then $M = G/P$, the quotient of a
semisimple algebraic group $G$ by a parabolic subgroup $P \subset
G$.
\end{conj}

Thus all the symplectic resolutions of the form $T^*M$ are
conjecturally covered by Example~\ref{homo.exa}. We note that if in
assumptions of Conjecture~\ref{dema} one in addition requires that
$Y$ has an isolated singularity, then one can prove that $M$ must be
a projective space -- this is S. Mori's famous theorem on smooth
varieties with ample tangent bundle.

Another natural source of symplectic resolution comes from global
holomorphic symplectic geometry. Given a projective holomorphic
symplectic manifold $\overline{X}$, one can sometimes construct a
projective birational contraction $\overline{X} \to \overline{Y}$;
the preimage $X \subset \overline{X}$ of any open affine $Y \subset
\overline{Y}$ is then an example of a symplectic resolution
according to Definition~\ref{reso.defn}. However, so far, all the
singularities obtained in this way are also covered by
Example~\ref{homo.exa} (prominent examples are the Mukai
Contraction, where $X=T^*\Pp^n$, and the O'Grady singularity, where
$X = T^*L$, $L$ the Grassmanian of Lagrangian $2$-planes in the
$4$-dimensional symplectic vector space).

\medskip

Let us now list some of the general properties of symplectic
singularities and symplectic resolutions. First of all, assume given
a symplectic singularity $Y$. The first observation is the
following.

\begin{lemma}[\cite{beau}]\label{ratio}
A symplectic singularity $Y$ is
necessarily canonical and rational.\endproof
\end{lemma}

\begin{corr}\label{no.h1}
For any fiber $F$ of a smooth resolution $\pi:X \to Y$ of a
symplectic singularity $Y$, the cohomology group $H^1(F_{an},\Z)$ is
trivial.
\end{corr}

\proof{} By proper base change, it suffices to prove that
$R^1\pi_*\Z_X=0$; this immediately follows from $R^1\pi_*\calo_X=0$
by considering the exponential exact sequence.
\endproof

To proceed, it turns out to be very productive to notice that $Y$
carries a natural structure of a Poisson scheme, \cite[Definition
1.2]{k1}. Namely, for every two local functions $f,g \in \calo_Y$ on
$Y$, we define
$$
\{f,g\} = \Theta \cntrct df \wedge dy
$$
on the smooth part $Y^{sm} \subset Y$, where $\Theta \in
H^0(Y^{sm},\Lambda^2\T)$ is the bivector on $Y^{sm}$ dual to the
symplectic form, and we note that since $Y$ is normal, any function
extends uniquely from $Y^{sm}$ to the whole $Y$. The advantage of
the Poisson bracket over the symplectic form is that the bracket is
perfectly well-defined on the singular locus of $Y$. This allows to
prove the following.

\begin{theorem}[{{\cite[Theorem 2.3]{k1}}}]\label{sympl.singu.th}
Every symplectic singularity $Y$ admits a finite stratification by
locally closed Poisson subschemes $Y_i \subset Y$ such that every
$Y_i$ is in fact smooth, and the induced Poisson structure on $Y_i$
comes from a symplectic form. All the closures $\overline{Y_i}$ are
also symplectic singularities. Moreover, for any closed point $y \in
Y$, the formal completion $\wh{Y}_y$ admits a decomposition
\begin{equation}\label{wein.eq}
\wh{Y}_y = Y_{y,0} \whtimes \wh{Y_i}_y
\end{equation}
into the Poisson scheme product of the formal germ of the stratum
$Y_i$ containing $y$ and a certain symplectic singularity $Y_{y,0}$.
\end{theorem}

The decomposition \eqref{wein.eq} into the product of a symplectic
stratum and a transversal slice is an algebraic version of the {\em
Weinstein decomposition} known in the $C^\infty$ Poisson geometry,
see \cite{wein}. Unfortunately, in algebraic geometry it only exists
after passing to the formal completion. To remedy the situation, it
would be very convenient to have one additional piece of structure
on a symplectic singularity -- namely, a $\gm$-action which dilates
the symplectic form.

\begin{defn}
An action of the algebraic group $\gm$ on a symplectic singularity
$Y$ is called {\em dilating} if it preserves the line $K \cdot
\Omega \in H^0(Y^{sm},\Omega^2_Y)$ and acts on this line via a
representation of weight $l > 0$ (in other words, $\lambda \cdot
\Omega = \lambda^l\Omega$ for some fixed integer $l > 0$ any
$\lambda \in \gm(K) = K^*$). A $\gm$-action on $Y$ is said to be
{\em positive-weight} if every finite-dimensional $\gm$-equivariant
subquotient of $H^0(Y,\calo_Y)$ decomposes into a sum of
representations of non-negative weights.
\end{defn}

\begin{conj}\label{gm.conj}
Every transversal slice $Y_{y,0}$ in \eqref{wein.eq} admits a
dilating po\-si\-ti\-ve\--weight action of $\gm$ such that $y \in
Y_{y,0}$ is its only fixed point.
\end{conj}

This conjecture is not as wild as one might suppose, since it is
actually possible to prove that $Y_{y,0}$ does admit a dilating
$\gm$-action whose only fixed point is $y$ (\cite[Theorem
2.4]{k1}). It is {\em not} true, however, that such an action always
has positive weights. It often happens in examples that $Y$ admits a
commuting $\gm$ action which is Hamiltonian, and even if one starts
with a positive-weight dilating action, composing it with a
commuting Hamiltonian action may produce a dilating $\gm$ action
whose weights are no longer positive. Effectively,
Conjecture~\ref{gm.conj} says that this is the only source of
problems: whenever a dilating action provided by \cite[Theorem
2.4]{k1} is not positive-weight, there is a commuting Hamiltonian
action which can be used to correct the weights.

In examples, a good $\gm$-action is always present. In the quotient
singularity case $Y=V/G$, the strata $Y_i$ and the transversal
slices $Y_{y,0}$ are also symplectic quotient singularities, and the
$\gm$-action on $V$ with weight $1$ induces a good dilating
action. In the nilpotent orbit case, the strata are smaller
nilpotent orbits, and it is possible to obtain a dilating
positive-weight $\gm$-action on transversal slice by an explicit
construction. The same is true in the quiver variety case.

The existence of a positive-weight dilating $\gm$-action is in fact
necessary for some of the results that we are going to describe;
since we cannot at present prove Conjecture~\ref{gm.conj}, it has to
be introduced as an additional hypothesis in the statements.


\medskip

Assume now given a symplectic resolution $X$ in the sense of
Definition~\ref{reso.defn}. Then the affine scheme $Y = \Spec
H^0(X,\calo_X)$ is a symplectic singularity, so that
Theorem~\ref{sympl.singu.th} applies. However, one can also prove
some stronger results concerning the topology of $X$.

\begin{theorem}[{{\cite[Corollary 2.8, Lemma 2.11, Theorem 2.12]{k1}}}]
  \label{topo}
Let $X$ be a symplectic resolution of $Y=\Spec H^0(X,\calo_X)$. Then
$X$ is semismall over $Y$ -- that is, $\dim X \times_Y X = \dim
X$. More exactly, for any stratum $Y_i \subset Y$ of codimension
$\codim Y_i = 2l$, its preimage $X_i \subset X$ has codimension
$\geq l$. Moreover, we have
$$
H^p(X,\Omega^q_X)=0
$$
whenever $p > q$, and for any fiber $E \in X$ of the map $X \to Y$,
the odd cohomology groups $H^{2p+1}(E_{an},\C)$ of the corresponding
complex-analytic space $E_{an}$ are trivial, while the even
cohomology groups $H^{2p}(E_{am},\C)$ carry a pure $\R$-Hodge
structure of weight $2p$ and type $(p,p)$. Finally, the symplectic
form $\Omega$ is exact in the formal neighborhood of any fiber $E
\in X$.
\end{theorem}

\proof[Sketch of the proof.] For any integer $i \geq 0$, denote by
$F^iH^\hdot_{DR}(X)$ the image of the cohomology
$H^\hdot(X,F^i\Omega^\hdot_X)$ with coefficients in the $i$-th term
of the stupid filtration on the de Rham complex
$\Omega^\hdot_X$. This is not quite the Hodge filtration, since $X$
is not compact, but it is functorial, and restricts to Hodge
filtration on compact fibers. By definition, the de Rham cohomology
class $[\Omega] \in H^2_{DR}(X)$ of the symplectic form $\Omega$
lies in $F^2H^2_{DR}(X)$. Since $Y$ has rational singularities, we
have $H^i(X,\calo_X)=0$ for $i \geq 1$, so that $H^\hdot_{DR}(X) =
F^1H^\hdot_{DR}(X)$, and the same is then true for the
complex-conjugate to the filtration $F^\hdot$. In particular,
$[\Omega] \in \overline{F^1H^2_{DR}(X)}$. By Hodge theory, this
implies that $[\Omega]$, hence also $\Omega$ itself restricts to $0$
on any fiber $E$.

Since the de Rham cohomology $H^\hdot_{DR}(\wt{E})$ of the formal
neighborhood $\wt{E}$ of the fiber $E \in X$ is isomorphic to the
cohomology of $E$ itself, $[\Omega]$ is also trivial in
$H^2_{DR}(\wt{E})$. Noting that $H^1(\wt{E},\calo_{\wt{E}})=0$ by
Lemma~\ref{ratio}, we deduce that $\Omega$ is exact on $\wt{E}$.

Let now $Y_i \subset Y$ be some smooth stratum, let $y \in Y_i$ be a
closed point, and let $E \subset X$ be its preimage. Since
$\Omega=0$ on $E$, for any tangent vector $\xi \in T_yY_i$ we have a
well-defined $1$-form $\alpha = \xi \cntrct \Omega$ on the smooth
part of $E$. By careful analysis of the construction of the Hodge
structure on the cohomology of $E$ -- this is a delicate point, for
details we refer the reader to \cite[Lemma 2.9]{k1} -- one checks
that $\alpha$ actually has a well-defined class in
$H^1(E_{an},\C)$. But this group vanishes by
Corollary~\ref{no.h1}. Again by Hodge theory, this implies that
$\alpha=0$ on the smooth part of $E$.

We conclude that the restriction of $\Omega$ to the smooth part of
the preimage $X_i$ is obtained by pullback from a symplectic form on
$Y_i$ (more careful analysis actually shows that this is the
symplectic from induced by the Poisson structure on $Y$). Since
$\Omega$ is non-degenerate on $X$, this gives the dimension estimates
$2\codim X_i \geq \codim Y_i$, $\dim X \times_Y X = \dim X$.

Now one can apply one of the vanishing theorems of \cite{ev}, see
\cite[Lemma 2.10]{k1}, to conclude that $H^p(X,\Omega^q_X)=0$
whenever $p+q \geq \dim X$, or, equivalently, whenever $p \geq q$
(since $X$ is symplectic, we have $\Omega^q_X \cong \Omega^{\dim
X-q}_X$). Again analyzing the Hodge structure on the fiber $E$, we
deduce the cohomological purity claim.  \endproof

\section{Deformations and quantization.}

As before, let us fix a symplectic resolution $X$ over a field $K$
of characteristic $0$. Next we would like to discuss deformation
theory of $X$.

Since $X$ is not compact, it usually does not have a reasonable
deformation theory as an algebraic variety (in particular, the
deformation space is infinite-dimensional). What one has to do is to
consider deformations of the pair $\langle X,\Omega\rangle$ -- that
is, deformations of $X$ together with the symplectic form
$\Omega$. For this deformation problem, the order-$1$ deformations
are controlled by the second de Rham cohomology group $H^2_{DR}(X)$
which is finite-dimensional. Moreover, in higher orders we have a
complete analog of the Bogomolov-Tian-Todorov Unobstructedness
Theorem.

\begin{theorem}[{{\cite[Theorem 1.1]{kv}}}]\label{defo}
The pair $\langle X,\Omega \rangle$ admits a universal formal
deformation $\langle \X,\Omega_{\X}\rangle/S$ whose base is the
completion of the vector space $H^2_{DR}(X)$ at the point $[\Omega]
\in H^2_{DR}(X)$.\endproof
\end{theorem}

In fact, for any deformation $\langle \X',\Omega_{\X'}\rangle/S'$,
the classifying map $S' \to S$ is the period map for the symplectic
form $\Omega$: it sends a point $s \in S'$ to $[\Omega_s] \in
H^2_{DR}(X)$, where $\Omega_s$ is the symplectic form on the fiber
$\X'_s$ over $s \in S'$, and we identify $H^2_{DR}(\X'_s) \cong
H^2_{DR}(X)$ by using the Gauss-Manin connection. Of course, the
deformation $\X/S$ is only formal, so that $s \in S'$ should be
understood her as an $A$-valued point for some Artin local
$K$-algebra $A$. The proof imitates the usual proof of the
Bogomolov-Tian-Todorov Lemma using the $T_1$-lifting principle of
Z. Ran (\cite{ran}, \cite{kawaran}).

\medskip

However, it turns out that it also makes sense to consider the {\em
non-commutative} deformations of $X$. Namely, we introduce the
following.

\begin{defn}\label{qua.defn}
A {\em quantization} of a Poisson variety $X$ is a sheaf $\calo_h$
of flat $K[[h]]$-algebras on $X$, complete in $h$-adic topology and
equipped with an isomorphism $\calo_h/h \cong \calo_X$, $f \mapsto
\overline{f}$, so that for any two local sections $f,g \in \calo_h$
we have
$$
fg-gf = h\{\overline{f},\overline{g}\} \mod h^2,
$$
where $\{-,-\}$ is the Poisson bracket on $X$.
\end{defn}

This can be applied to our symplectic resolution $X$, or to some
deformation $\X'/S'$ with the Poisson structure induced by the
relative symplectic form (functions in $\calo_{S'}$ are in the
Poisson center). This allows one to give a complete classification
of quantization.

\begin{theorem}[{{\cite[Theorem 1.8, Lemma 6.4]{bk1}}}]\label{qua}
Assume given a symplectic resolution $X$, and let $\X/S$ be its
universal symplectic deformation provided by
Theorem~\ref{defo}. Then there exists a canonical quantization
$\wt{\calo}_h$ of the Poisson variety $\X$ which is universal in the
following sense: for any quantization $\calo_h$ of the symplectic
resolution $X$, there exists a unique section
$$
s:\Spec K[[h]] \to \Spec K[[h]] \whtimes S
$$
of the projection $\Spec K[[h]] \whtimes S \to \Spec k[[h]]$ such
that $s^*(\wt{\calo}_h)$ is isomorphic to $\calo_h$.
\end{theorem}

Algebraically speaking, quantizations of $X$ exist, and they are
classified up to an isomorphism by power series $H^2_{DR}(X)[[h]]$
with coefficients in $H^2_{DR}(X)$ and leading term $[\Omega]$.

We would like to explain briefly why Theorem~\ref{qua} holds. For
this, we need to start with the local theory.

\subsection{Local theory.}
Let $A=K[[x_1,\dots,x_n,y_1,\dots,y_n]]$ be the algebra of functions
on the formal neighborhood of $0$ in the $2n$-dimensional vector
space. Equip $\Spec A$ with the standard symplectic form $\Omega=dx_1
\wedge dy_1 + \dots + dx_n \wedge dy_n$ and consider the associated
Poisson structure. Then $A$ has a standard quantization $D$ given by
\begin{equation}\label{D.eq}
D=K[[x_1,\dots,x_n,y_1,\dots,y_n,h]]/\{x_ix_j-x_jx_i, y_iy_j-y_jy_i,
x_iy_j-y_jx_i-\delta_{ij}h\},
\end{equation}
where $\delta_{ij}$ is the Kronecker delta-symbol. Denote by
$\Aut(D)$ the algebraic group of $K[[h]]$-linear automorphisms of
the algebra $D$. For any $l \geq 1$, denote by $\Aut_{\geq l}(D)
\subset \Aut(D)$ the subgroup of automorphisms which are equal to
identity on $D/h^lD$.

\begin{prop}\label{qua.loc}
\begin{enumerate}
\item For any other symplectic form $\Omega'$ on $\Spec A$, the pair
$\langle A,\Omega \rangle$ is isomorphic to $\langle A, \Omega'
\rangle$.
\item Any quantization of $\langle A,\Omega \rangle$ is
isomorphic to $D$. 
\item The subgroup $\Aut_{\geq 1}(D) \subset \Aut(D)$ coincides with
the subgroup of inner automorphisms of the algebra $D$, and we have a
central group extension
\begin{equation}\label{centra}
\begin{CD}
1 @>>> K[[h]]^* @>>> D^* @>>> \Aut_{\geq 1}(D) @>>> 1.
\end{CD}
\end{equation}
\end{enumerate}
\end{prop}

\proof[Sketch of a proof.] \thetag{i} is the Darboux Theorem; all
the standard proofs work in our formal algebraic situation. For
\thetag{ii} and \thetag{iii}, we first prove the following.

\begin{lemma}\label{hoch}
For any quantization $D'$ of the Poisson algebra $A$, the relative
Hochschild cohomology groups $HH^l_{K[[h]]}(D')$ of the algebra $D'$
over $K[[h]]$ are annihilated by $h$ for $l \geq 1$, while
$HH^0_{K[[h]]}(D') \cong K[[h]]$.
\end{lemma}

\proof[Sketch of a proof.] Consider the spectral sequence computing
$HH^\hdot_{K[[h]]}(D')$ associated to the $h$-adic filtration on
$D'$. Its first term $E_1^{\hdot,\hdot}$ is
$HH^\hdot_{K[[h]]}(\gr^\hdot D')$, and the associated graded
quotient $\gr^\hdot D'$ with respect to the $h$-adic filtration by
definition coincides with $A[[h]]$. Therefore we have
$$
E^{\hdot,\hdot}_1 \cong \Lambda^\hdot(\T(A))[[h]],
$$
the algebra of polyvector fields on $\Spec A$. One checks easily
that the differential in the spectral sequence is given by 
$$
d_1(\alpha) = h[\Theta,\alpha],
$$
where $\Theta \in \Lambda^2(\T(A))$ is the Poisson bivector, and
$[-,-]$ is the Schouten bracket of polyvector fields. If we use the
symplectic form $\Omega$ on $A$ to identify $\T(A) \cong
\Omega^1(A)$, then $\Lambda^\hdot\T(A)$ becomes identified with the
algebra $\Omega^\hdot(A)$ of differential forms, and the
differential $d_1$ becomes identified with $hd$, where
$d:\Omega^\hdot(A) \to \Omega^{\hdot+1}(A)$ is the de Rham
differential. But by the Poincar\'e Lemma, $\Spec A$ has no higher
de Rham cohomology, and $H^0_{DR}(A) \cong K$. Therefore for
dimension reasons, the spectral sequence degenerates already at
$E^{\hdot,\hdot}_2$, and the term $E^{\hdot,\hdot}_2 \cong
E^{\hdot,\hdot}_\infty$ is as required by the statement.
\endproof

Now, since all the groups in Proposition~\ref{qua.loc}~\thetag{iii}
are unipotent, it suffices to prove the corresponding statement for
Lie algebras, where it immediately follows from Lemma~\ref{hoch}:
since $HH^1(D)$ is annihilated by $h$, every derivation of the
algebra $D$ which is divisible by $h$ is inner, and the center of
the algebra $D$ is $HH^0(D) = K[[h]]$. To prove \thetag{ii}, we take
a different quantization $D'$, and we prove by induction that
$D/h^lD \cong D'/h^lD'$ for any $l$. Indeed, $D/h^lD$ is an
extension of the algebra $D/h^{l-1}D$ by $A/K$, and as as such, it
is given by an extension class $\beta$ in
$$
HH^2_{K[h]/h^{l-1}}(D/h^{l-1}D,A/K) \cong \Omega^2(A).
$$
If $l=2$, then this class $\beta$ is exactly the Poisson bivector
for $A$, and therefore it is the same for $D$ and for $D'$. If $l
\geq 2$, then $\beta$ comes from a class in
$HH^2(D/h^{l-1}D,D/h^2D)$, so that it survives in the $E_2$-term of
the $h$-adic spectral sequence. By the proof of Lemma~\ref{hoch},
this means that the corresponding form $\beta \in \Omega^2(A)$ must
be closed. The same is true for the class $\beta' \in \Omega^2(A)$
associated to $D'/h^lD'$. But by the Poincar\'e Lemma, all closed
forms on $A$ are exact, $\beta = \beta' + d\alpha$; therefore one
can change an isomorphism $D/h^{l-1}D \cong D'/h^{l-1}D'$ by
composing it with $\exp(1 + h^{l-1}\alpha)$ so that we have
$\beta=\beta'$ and $D'/h^lD' \cong D/h^lD$.  \endproof

For any $l \geq 0$, denote the quotient $\Aut(D)/\Aut_{\geq l+1}(D)$
by $\Aut^l(D)$. By Proposition~\ref{qua.loc}~\thetag{ii}, every
symplectic automorphism of $A$ lifts to an automorphism of $D$, so
that $\Aut^0(D)$ is the group of symplectic automorphisms of
$A$. For any $l \geq 1$, we have a group extension
\begin{equation}\label{et.n}
\begin{CD}
1 @>>> A/K @>>> \Aut^l(D) @>>> \Aut^{l-1}(D) @>>> 1,
\end{CD}
\end{equation}
where $A/K$ on the left-hand side is considered as the additive
group. By definition, $\Aut^l(D)$ acts on $D/h^{l+1}D$, and we have
an injective map
$$
\Aut^l(D) \hookrightarrow \Aut(D/h^{l+1}D).
$$
However, this map is not surjective -- not every $K[[h]]$-linear
automorphism of $D/h^{l+1}D$ lifts to an automorphism of $D$. To
describe the groups $\Aut^l(D)$ exactly, we introduce the following.

\begin{defn}\label{qua.alg.defn}
A {\em quantized algebra} over the field $K$ is an associative
algebra $B$ over $K[[h]]$ equipped with a Lie bracket $\{-,-\}$ such
that $\{a,-\}$ is a derivation of $B$ for any $b \in B$, and for any
$a,b \in B$ we have
$$
ab-ba = h\{a,b\}.
$$
\end{defn}

If a quantized algebra $B$ is flat over $K[[h]]$, then the bracket
can be recovered uniquely from the product, and
Definition~\ref{qua.alg.defn} simply requires that the commutator
$ab-ba$ is divisible by $h$ (in other words, $B/hB$ is
commutative). Conversely, a quantized algebra annihilated by $h$ is
the same thing as a Poisson algebra. In general,
Definition~\ref{qua.alg.defn} interpolates between the two
notions.

\begin{lemma}\label{qua.alg.le}
For any $l \geq 0$, the group $\Aut^l(D)$ is the group of
$K[[h]]$-linear automorphisms of the quantized algebra $D/h^{l+1}D$.
\end{lemma}

\proof{} Left to the reader. \endproof

\subsection{Globalization by formal geometry.} Assume now given a
smooth symplectic variety $X$ of dimension $2n$ over $K$. To obtain
global results about quantization of $X$, we use the technique
called {\em formal geometry} (\cite{GeKa}). Namely, we note that $X$
defines a completely canonical {\em variety of formal coordinate
systems} $\M_{coord}$ whose points are pairs $\langle x,\phi
\rangle$ of a point $x \in X$ and an isomorphism
$\phi:\wh{\calo_{X,x}} \cong A$ between the standard algebra $A$ and
the algebra of functions on the formal neighborhood of $x$ in
$X$. Forgetting $\phi$ gives a map
$$
\M_{coord} \to X,
$$
and one check easily that $\M_{coord}$ is a torsor over $X$ with
respect to the group $\Aut(A)$ of all continuous automorphisms of
the algebra $A$. The symplectic structure on $X$ defines a subtorsor
$$
\M_0 \subset \M_{coord}
$$
of {\em symplectic} coordinate systems -- that is, such pairs
$\langle x,\phi \rangle$ that $\phi$ is compatible with the
symplectic form. By Proposition~\ref{qua.loc}~\thetag{i}, the
forgetful map $\M_0 \to X$ surjective, and $\M_0$ is a torsor with
respect to the group $\Aut^0(D) \subset \Aut(A)$ of all symplectic
automorphisms of $A$ -- in other words, we have a restriction of the
torsor $\M_{coord}$ to the subgroup $\Aut^0(D) \subset \Aut(A)$.

Analogously, every quantization $\calo_h$ of the symplectic manifold
$X$ defines the space $\wt{\M}$ of pairs $\langle
x,\wt{\phi}\rangle$, $x \in X$, $\wt{\phi}$ is an isomorphism
between the natural completion $\wh{\calo_h}$ at $x \in X$ and the
standard quantization $D$. By Proposition~\ref{qua.loc}~\thetag{ii},
$\wt{\M}$ is a restriction (or maybe, a lifting would be a better
term) of the $\Aut^0(D)$-torsor $\M_0$ to the group $\Aut(D)$ with
respect to the natural map $\Aut(D) \to \Aut^0(D)$.

\medskip

Thus to any quantization $\calo_h$, we associate a restriction of
the $\Aut^0(D)$-torsor $\M_0$ to the group $\Aut(D)$.

\medskip

This is not a one-to-one correspondence, but it can be made into one
if we enlarge the automorphism groups. Namely, every automorphism of
the algebra $D$ must necessarily preserve the maximal ideal $\m
\subset D$ generated by $x_1,\dots,x_n,y_1,\dots,y_n,h \in D$. Thus
the Lie algebra of the group $\Aut(D)$ is the algebra $\Der_0(D)$ of
all $K[[h]]$-linear derivations $D \to D$ which preserve $\m \subset
D$. It lies naturally within a larger algebra $\Der(D)$ of {\em all}
$K[[h]]$-linear derivations, but those derivations which do not
preserve $\m$ cannot be integrated to actual automorphisms. We have
a similar picture for $\Aut(A)$. To take account of the
non-integrable derivations, one has to either consider groups which
are partially formal, or, which is simpler, to work with the
so-called {\em Harish-Chandra pairs} $\langle G,\h \rangle$ of an
algebraic group $G$, a Lie algebra $\h$ equipped with an action of
$G$, and a $G$-equivariant Lie algebra embedding $\g \hookrightarrow
\h$, where $\g$ is the Lie algebra of the group $G$. When the
appropriate notion of a $\langle G,\h \rangle$-torsor is introduced
(see e.g. \cite[Definition 2.3]{bk1}, or \cite[Section 2.6]{BeDr}),
$\M_{coord}$ becomes a $\langle \Aut A,\Der A\rangle$-torsor over
$X$, $\M_0$ becomes a $\langle \Aut^0(D),\Der^0(D)\rangle$-torsor,
and we have the following.

\begin{prop}[{{\cite[Lemma 3.4]{bk1}}}]\label{tors}
Let $X$ be a smooth symplectic variety over $K$ of dimension
$2n$. Then there exists a natural one-to-one correspondence between
the isomorphism classes of quantizations $\calo_h$ of the variety
$X$ and the isomorphism classes of liftings $\wt{\M}$ of the
symplectic coordinate torsor $\M_0$ from the Harish-Chandra pair
$\langle \Aut^0(D),\Der^0(D)\rangle$ to the Harish-Chandra pair
$\langle \Aut(D),\Der(D)\rangle$.\endproof
\end{prop}

An analogous statement holds in the relative situation -- for a
variety $\X$ smooth and symplectic over a base $S$ and of relative
dimension $2n$.

\medskip

With the use of Proposition~\ref{tors}, the problem of finding and
classifying quantizations reduces to the study of torsors. This can
be done step-by-step using the extensions \eqref{et.n}. By the
standard obstruction theory, at each step, the obstruction to
lifting an $\Aut^l(D)$-torsor to $\Aut^{l+1}(D)$ lies in the group
$$
H^2(X,\E),
$$ 
where $\E$ is the vector bundle on $X$ associated to the given
$\Aut^l(D)$-torsor via the action of $\Aut^l(D)$ on $A/K$. This
action factors through $\Aut^0(D)$; therefore $\E$ does not really
depend on the $\Aut^l(D)$-torsor. An easy computation shows that $\E
\cong J^\infty\calo_X/\calo_X$, the quotient of the jet bundle
$J^\infty\calo_X$ of the structure sheaf by the structure sheaf
itself. When we pass to Harish-Chandra pairs, the cohomology gets
replaced with de Rham cohomology (see \cite{bk1}; note that both
$\calo_X$ and $J^\infty\calo_X$ carry natural flat
connections). Thus the group that contains the obstruction fits as
the middle term into an exact sequence
$$
\begin{CD}
H^2_{DR}(X,J^\infty\calo_X) @>>> H^2_{DR}(X,J^\infty\calo_X/\calo_X)
@>>> H^3_{DR}(X)
\end{CD}
$$
Moreover, since the central extension \eqref{centra} is obviously
compatible with the filtration $\Aut_{\geq \idot}(D)$, one shows
that the obstruction actually comes from an element in
$H^2_{DR}(X,J^\infty\calo_X)$. But by the universal property of the
jet bundle, this group coincides with $H^2(X,\calo_X)$, which is
trivial for a symplectic resolution $X$. Thus there are no
obstructions. Analogously, we see that possible liftings at each
step are classified by elements of the group
$$
H^1_{DR}(X,J^\infty\calo_X/\calo_X),
$$
and this is isomorphic to $H^2_{DR}(X)$, as required in
Theorem~\ref{qua}.

\begin{remark}
In fact, one can use a parallel argument to study commutative
symplectic deformations and prove Theorem~\ref{defo}. To do this,
one replaces the standard quantization $D$ by the standard (trivial)
symplectic deformation $A[[h]]$, and considers its automorphism
group. An exact sequence completely analogous to \eqref{et.n}
exists, where $A/K$ on the left-hand side consists of exponentials
of Hamiltonian vector fields on $\Spec A$.
\end{remark}

\section{Positive characteristic case.}

Interestingly, the quantization theory for symplectic resolutions
can be developed still further and in a slightly unexpected
direction. Namely, assume now that everything is defined over a
perfect base field $k$ of {\em positive} characteristic $\cchar k =
p > 0$.

\subsection{New phenomena.}
Our definition of a quantization, Definition~\ref{qua.defn}, makes
sense in positive characteristic without any changes, and so does
the standard example of a quantization, the algebra $D$ of
\eqref{D.eq}. However, in $\cchar p$ we immediately see a new
feature: the algebra $D$ has a large center. The $p$-th powers
$x_i^p,y_j^p$ of the generators become central in $D$. This
observation motivates the following.

\begin{defn}[{{\cite[Definition 1.4]{bk3}}}]\label{pqua.defn}
A {\em Frobenius-constant} quantization of a Poisson scheme $X$ over
a field $k$ of characteristic $p$ is a pair of a quantization
$\calo_h$ of $X$ in the sense of Definition~\ref{qua.defn} and a map
$s:\calo_X \to \calo_h$ such that for any $f \in \calo_X$, $s(f)$ is
central in $\calo_h$ and satisfies
$$
s(f) = f^p \mod h^{p-1}.
$$
\end{defn}

In other words, the natural surjection $\calo_h \to \calo_X$ should
be split by the map $s$ on the subalgebra $\calo_X^p \subset
\calo_X$ of $p$-th powers of functions on $X$.

\begin{remark}
Definition~\ref{pqua.defn} first appeared in \cite{bk2} with a
weaker condition, $s(f)=f^p \mod h$, not $\mod h^{p-1}$. The need
for a stronger condition was realized in \cite{bk3}.
\end{remark}

One would like to prove a version of Theorem~\ref{qua} for
Frobenius-constant quantizations in positive characteristic. We can
start in the same way, and we immediately notice that the situation
is somewhat simpler: instead of considering the whole algebra $A$,
it is enough to consider the truncated algebra
\begin{equation}\label{barA.eq}
\barA = A/\{x_i^p,y_i^p\}
\end{equation}
and its standard quantization
$$
\barD = D/\{x_i^p,y_i^p\}.
$$
Again, we consider the algebraic group $\Aut(\barD)$, its subgroups
$\Aut_{\geq l}(\barD)$, and the quotient groups $\Aut^l(\barD)$. But
since $D$ is a finitely-generated $k[[h]]$-module, these groups are
now finite-dimensional. For the same reason, there is no need to
consider Harish-Chandra pairs: the non-integrable automorphisms can
be included into the group $\Aut(\barD)$ itself as a certain
non-reduced part (formally, $\Aut(\barD)$ is the group scheme over
$k$ which represents the functor $T \mapsto \Aut_{T[[h]]}(D \otimes
\calo_T)$).

\medskip

However, there is a price to pay for this simplification, and it is
Proposition~\ref{qua.loc}. Out of its three claims, only
\thetag{iii} survives in positive characteristic essentially in the
same form.

\begin{lemma}\label{azu}
For any Frobenius-constant quantization $\barD'$ of the Poisson
algebra $\barA$, the tensor product $\barD' \otimes_{k[[h]]} k((h))$
is a matrix algebra over the Laurent power series field $k((h))$,
and we have a group extension
\begin{equation}\label{centra.p}
\begin{CD}
1 @>>> k[[h]]^* @>>> \barD^{'*} @>>> \Aut_{\geq 1}(\barD') @>>> 1,
\end{CD}
\end{equation}
where $\Aut_{\geq 1}(\barD')$ is the algebraic group of all
$k[[h]]$-linear automorphisms of the algebra $\barD'$ which are
equal to identity on $A = \barD'/h$.
\end{lemma}

\proof[Sketch of a proof.] The first claim is \cite[Lemma
3.2]{bk3}; the second claim immediately follows from the
Skolem-Noether Theorem.
\endproof

\begin{remark}
Lemma~\ref{azu} is true even for a weaker notion of a
Frobenius-constant quantization used in \cite{bk2}.
\end{remark}

Proposition~\ref{qua.loc}~\thetag{ii} also has a
positive-characteristic counterpart -- if one uses the strong notion
of a Frobe\-ni\-us\--constant quantization, and imposes an
additional technical assumption (Definition~\ref{regu}). However, the
proof is much more delicate, since Lemma~\ref{hoch} no longer holds.

Indeed, our proof of Lemma~\ref{hoch} goes through in positive
characteristic without any changes up to the point where one needs
to invoke the Poincar\'e Lemma. But the Poincar\'e Lemma now is
false: the higher de Rham cohomology groups $H^\hdot_{DR}(A)$ are
not trivial. Conversely, they are very large -- there exists the
so-called {\em Cartier isomorphism} $C:H^\hdot_{DR}(A) \cong
\Omega^\hdot(A^p)$ which identifies the de Rham cohomology groups of
$\Spec A$ with the spaces of differential forms on $\Spec A^p$
(recall that the subring $A^p \subset A$ of $p$-th powers of
elements of $A$ is isomorphic to $A$ as an abstract ring). For the
reduced algebra $\barA$, the Cartier map identifies de Rham
cohomology algebra $H^\hdot_{DR}(\Spec \barA)$ with the exterior
algebra $\Lambda^\hdot W$, where $W= T^*_o\Spec A$ is the cotangent
vector space at the closed point $o \in \Spec A$.

For the same reason, Proposition~\ref{qua.loc}~\thetag{i} -- in
other words, the Darboux Theorem -- completely breaks down in
positive characteristic. For example, while the standard symplectic
form $\Omega$ on $A$ satisfies $C(\Omega)=0$, it is very easy to
construct a symplectic form $\Omega'$ with a $C(\Omega')$ not equal
to $0$, or in fact equal to an arbitrary prescribed non-zero
$2$-form on $A^p$. Since the Cartier map is canonical, $A$ equipped
with such a form $\Omega'$ cannot be isomorphic to $\langle
A,\Omega\rangle$.

\begin{remark}
Since Hochschild cohomology is Morita-invariant, a matrix algebra
has the same Hochschild cohomology as its center; therefore
Lemma~\ref{azu} implies that $HH^\hdot_{k[[h]]}(\barD')$ is
annihilated by some power of $h$, just as in Lemma~\ref{hoch}. What
happens is that there is a second non-trivial differential in the
spectral sequence, namely, the differential at the term $E^p$. After
that term, the spectral sequence does degenerate, and we still get
the statement of Lemma~\ref{hoch}, but the higher Hochschild
cohomology groups are only annihilated by $h^{p-1}$, not by
$h$. This, in particular, explains why, if we want
Proposition~\ref{qua.loc} to hold, we have to introduce a stronger
assumption on $s(f)$ in Definition~\ref{pqua.defn}.
\end{remark}

\subsection{Restricted structures.}
To recover the full statement of Proposition~\ref{qua.loc} in
positive characteristic, and then prove a version of
Proposition~\ref{tors}, we need to equip our quantizations with an
additional structure.

Denote by $Q(x,y)$ the free quantized algebra in the sense of
Definition~\ref{qua.alg.defn} generated by two elements $x$,
$y$. Explicitly, we have
$$
Q(x,y) = \bigoplus W^\hdot_{PBW}T^\hdot(x,y),
$$
where $T^\hdot(x,y)$ is the free associative algebra generated by
$x$ and $y$, and $W^\hdot_{PBW}$ is the Poincar\'e-Birkhoff-Witt
increasing filtration (see \cite[Subsection 1.2]{bk3} for
details). The quantization parameter $h$ acts on $Q(x,y)$ by the
natural embedding $W^\hdot_{PBW}T^\hdot(x,y) \to
W^{\hdot+1}_{PBW}T^\hdot(x,y)$. We call elements of the algebra
$Q(x,y)$ {\em quantized polynomials} in variables $x,y$.

\begin{lemma}
Assume that the base field $k$ is of characteristic $\cchar k = p >
0$. Then there exists quantized polynomials $F(x,y)$, $P(x,y)$ in
$x,y$ such that
$$
(x+y)^p - x^p - y^p = h^{p-1}F(x,y) \qquad (xy)^p -x^py^p =
h^{p-1}P(x,y).
$$
\end{lemma}

\proof{} Immediately follows from \cite[Lemma 1.3]{bk3}.
\endproof

\begin{defn}\label{restr.defn}
A {\em restricted quantized algebra} $A$ is a quantized algebra $A$
over a field $k$ of characteristic $\cchar k = p > 0$ equipped with
an additional operation $x \mapsto x^{[p]}$ such that $h^{[p]}=h$
and
\begin{equation}\label{restr.eq}
\begin{aligned}
\{x^{[p]},y\} &= (\ad x)^p(y),\\
(x+y)^{[p]} &= x^{[p]} + y^{[p]} + F(x,y),\\
(xy)^{[p]} &= x^py^{[p]} + x^{[p]}y^p - h^{p-1}x^{[p]}y^{[p]} +
P(x,y),
\end{aligned}
\end{equation}
where we denote by $\ad x:A \to A$ the endomorphism given by $y
\mapsto \{x,y\}$.
\end{defn}

\begin{remark}
The polynomial $F(x,y)$ is a well-known Lie polynomial in $x,y$
which can also be given by an explicit formula (Jacobson formula,
\cite[II, \S 7.3, D\'efinition 3.1]{dem2}). The first two equations
of \eqref{restr.eq} define the standard notion of a {\em restricted
Lie algebra}. The standard example of a restricted Lie algebra is
the algebra of vector fields on a scheme or, more generally, the
algebra of all derivations of an associative algebra $B/k$ (the
restricted power operation sends a derivation $D:B \to B$ to its
$p$-th power $D^p$, which, as one checks easily, is also a
derivation).
\end{remark}

Recall that a quantized algebra $A$ which is $h$-adically complete
and has no $h$-torsion is the same as a quantization in the sense of
Definition~\ref{qua.defn} (of the quotient $A/hA$). The notion of a
restricted quantized algebra plays the same role for
Frobenius-constant quantizations.

\begin{lemma}
A restricted quantized algebra $A$ which is $h$-adically complete
and has no $h$-torsion is the same as a Frobenius-constant
quantization of the quotient $A/hA$.
\end{lemma}

\proof[Sketch of a proof.] Given a restricted quantized algebra $A$,
we define the splitting map $s:A/hA \to A$ by
\begin{equation}\label{s.p}
s(a) = a^p - h^{p-1}a^{[p]}.
\end{equation}
The first equation of \eqref{restr.eq} guarantees that this is a
central map from $A$ to itself, and the last two equations insure
that it is an algebra map; since $h^{[p]}=h$, this map vanishes on
$hA \subset A$. Conversely, by our assumption, for any
Frobenius-constant quantization the difference $s(a)-a$ is divisible
by $h^{p-1}$ for any $a \in A$; since $A$ has no $h$-torsion,
\eqref{s.p} uniquely defines a restricted power operation $a \mapsto
a^{[p]}$, and the conditions on the map $s$ insure that
\eqref{restr.eq} is satisfied.
\endproof

On the other hand, if $h=0$ on a restricted algebra $A$, then $A$ is
a Poisson algebra. In this case Definition~\ref{restr.defn} gives a
notion of a {\em restricted Poisson algebra}. The first two
equations of \eqref{restr.eq} insure that $A$ equipped with the
Poisson bracket is a restricted Lie algebra in the usual sense, and
the last equation gives a compatibility condition between the
restricted power operation and the multiplication which, it seems,
was first introduced in \cite{bk3}. The main source of restricted
Poisson algebras is the following.

\begin{prop}[{{\cite[Theorems 1.11,1.12]{bk3}}}]\label{restr.prop}
Let $X = \Spec A$ be a smooth affine scheme over a field $k$ of
characteristic $\cchar k = p > 0$, and assume that $X$ is equipped
with a symplectic form $\Omega$. Then the following conditions are
equivalent.
\begin{enumerate}
\item We have $C([\Omega])=0$.
\item The form $\Omega$ is exact, $\Omega = d\alpha$ for some
  $\alpha \in \Omega^1_X$.
\item The Lie subalgebra $H \subset \T(A)$ of {\em Hamiltonian}
  vector fields on $X$ is closed with respect to the natural
  restricted power operation on $\T(A)$.
\item The Poisson algebra $A$ admits a restricted Poisson structure.
\end{enumerate}
Moreover, restricted Poisson structures in \thetag{iv} are in
one-to-one correspondence with $1$-forms $\alpha$ in \thetag{ii}
considered modulo exact $1$-forms, $\alpha \sim \alpha + df$ for any
$f \in A$.\endproof
\end{prop}

In particular, we see that the algebra $\barA$ equipped with the
standard symplectic form has a restricted Poisson structure.

\subsection{Quantization.}
We can now give the main results about quantization in positive
characteristic. First of all, we need the following technical notion
introduced in \cite{bk3}

\begin{defn}\label{qua.base.defn}
A {\em good quantization base} $B$ is a complete local $k$-algebra
$B$ with residue field $k$ equipped with an element $h$ in its
maximal ideal $\m \subset B$ and an additive operation $B \to B$, $b
\mapsto b^{[p]}$ such that $s:B \to B$ given by $s(b) = b^p -
h^{p-1}b^{[p]}$ is an algebra map.
\end{defn}

In other words, a good quantization base is a {\em commutative}
restricted quantized algebra, with a completeness condition. The
quotient $B/hB$ is then a complete local restricted Poisson algebra
whose Poisson bracket is trivial. The restricted power operation on
$B/hB$ need {\em not} be trivial (and in applications, it is
not). However, since $\{-,-\}=0$ tautologically, we must have
$b^{[p]}=K(b)$ for some additive map $K:B/h \to B/h$ which satisfies
$K(ab) = a^pK(b)+K(a)b^p$. In \cite{bk3}, such maps are called {\em
Frobenius-derivations}.

Given a quantization base $B$, by a {\em restricted quantized
algebra} $A$ over $B$ we will understand a quantized algebra over
$B$ equipped with a restricted structure in such a way that the
natural central embedding $B \to A$ is compatible with the
restricted structures.

We will also need the following notions.

\begin{defn}\label{regu}
A Frobenius-constant quantization $\langle \calo_h,s \rangle$ of a
scheme $X/k$ is called {\em regular} if for any local section $f \in
\calo_X$ such that $f^p=0$, we have $s(f)=0 \mod h^{p-1}$. A
restricted quantized algebra $A$ is called regular if for any $a \in
A$ with $a^p = 0 \mod h$, we have $a^{[p]} = 0 \mod h$.
\end{defn}

\begin{remark}
Regularity is a technical condition needed to study non-reduced
algebras such as the algebra $\barA$. We note that for reduced
algebras -- in particular, for algebras of functions on a smooth
algebraic variety -- this condition is tautologically satisfied.
\end{remark}

\begin{defn}
A {\em small Dieudonn\'e module} $I$ over $k$ is a $k$-vector space
equipped with an additive operation $I \to I$, $a \mapsto a^{[p]}$
which is Frobenius-semilinear, that is, $(\lambda a )^{[p]} =
\lambda^pa^{[p]}$ for any $\lambda \in k$, $a \in I$.
\end{defn}

Small Dieudonn\'e modules obviously form an abelian category. For
any good quantization base $B$ with an ideal $I \subset B$ such that
$\m_B \cdot I = 0$, the restricted power operation on $B$ induces a
structure of a small Dieudonn\'e module on $I$.

\begin{prop}[{{\cite[Proposition 3.8, Lemmas 3.10, 3.11]{bk3}}}]
  \label{restr.drb}
Assume given a good quantization base $B$ with maximal ideal $\m
\subset B$. Then there exists a unique regular restricted quantized
algebra $A^B$ over $B$ whose quotient $A^B/\m A^B$ is isomorphic as
a Poisson algebra to the standard Poisson algebra $\barA$ defined in
\eqref{barA.eq}. Moreover, for any restricted ideal $I \subset B$
such that $\m \cdot I = 0$ we have a natural extension of algebraic
groups
\begin{equation}\label{restr.ext}
\begin{CD}
1 @>>> \HH(\barA,I) @>>> \Aut(A^B) @>>> \Aut(A^{B/I}) @>>> 1,
\end{CD}
\end{equation}
where $\Aut(A^B)$ is the group of $B$-linear automorphisms of the
restricted quantized algebra $A^B$, $\Aut(A^{B/I})$ is the group of
$B/I$-linear automorphisms of the restricted quantized algebra
$A^{B/I}$, and $\HH(\barA,I)$ is a certain commutative algebraic
group which only depends on the small Dieudonn\'e module $I \subset
B$. Finally, if $I = \m$, then the group extension \eqref{restr.ext}
is a semi-direct product.\endproof
\end{prop}

This Proposition is the positive characteristic counterpart of
Proposition~\ref{qua.loc}, together with Lemma~\ref{qua.alg.le}. We
will not need the precise form of the group $\HH(\barA,I)$, see
\cite[Definition 1.16, Subsection 3.3]{bk3}. Let us just say that
$\HH(\barA,-)$ is an exact functor from small Dieudonn\'e modules to
commutative algebraic groups, and that the algebraic group
$\Aut(\barA)$ of all automorphisms of the algebra $\barA$ acts
naturally on $\HH(\barA,I)$ for any $I$ (the action of
$\Aut(A^{B/I})$ on $\HH(\barA,I)$ coming from the extension
\eqref{restr.ext} is induced by this action via the natural map
$\Aut(A^{B/I}) \to \Aut(\barA)$). Moreover, there are two particular
cases that we will need.
\begin{enumerate}
\item If the restricted structure on $I$ is trivial, $a^{[p]} = 0$
for any $a \in I$, then $\HH(\barA,I)$ is the $k$-vector space
$(\barA/k) \otimes_k I$ considered as an (additive) algebraic group.
\item If $I = k \cdot t$ for some element $t \in I$ with $t^{[p]} =
t$, then $\HH(\barA,I) = \barA^*/k^*$ (with the natural structure
of an algebraic group).
\end{enumerate}

To state the global quantization result, let us from now on, for any
scheme $X/k$, denote $X^\tw = \langle X,\calo_X^p \rangle$ -- that
is, $X$ with the subsheaf $\calo_X^p \subset \calo_X$ of $p$-th
powers as the structure sheaf. Let us denote by $\Fr:X \to X^\tw$
the natural map (if $X$ is reduced, $X^\tw$ is isomorphic to $X$,
and $\Fr$ becomes the Frobenius map). For any restricted Poisson
scheme $X/k$ and a good quantization base $B$ with maximal ideal $\m
\subset B$, by a {\em $B$-quantization} $\calo_B$ of $X$ we will
understand a sheaf of flat and complete restricted quantized
$B$-algebras $\calo_B$ on $X$ equipped with a restricted Poisson
isomorphism $\calo_B/\m_B \cong \calo_X$.

\begin{theorem}[{{\cite[Proposition 1.22]{bk3}}}]\label{pqua.thm}
Assume given a good quantization base $B$ in the sense of
Definition~\ref{qua.base.defn} and a smooth symplectic variety
$X/k$.  Assume also that $H^i(X,\calo_X)=0$ for $i=1,2,3$. Then the
isomorphism classes of $B$-quantizations of $X$ are in one-to-one
correspondence with elements of the \'etale cohomology group
\begin{equation}\label{etal}
H^1_{et}(X^\tw,\Loc(\HH(\barA,\m/\m^2))),
\end{equation}
where $\Loc(\HH(\barA,I))$ is an \'etale sheaf on $X^\tw$ which
only depends on the $\Aut(\barA)$-action on the algebraic group
$\HH(\barA,I)$. In particular, Frobenius-constant quantizations --
or, equivalently, $k[[h]]$-quantizations -- of $X/k$ are classified
up to an automorphism by elements of the group
$$
H^1_{et}(X^\tw,(\Fr_*\calo_X)^*/(\calo_{X^\tw})^*).
$$
Finally, every $B/h$-quantization of the scheme $X$ extends to a
$B$-quantization.
\end{theorem}

\proof[Sketch of the proof.] As in the proof of Theorem~\ref{qua},
we use the technique of formal geometry. To any smooth variety $X/k$
of dimension $2n$ one associates the torsor
$$
\M_{coord}(X) = \Maps(\Spec \barA,X)
$$
of \'etale maps $\Spec \barA \to X$ -- or, informally speaking, of
formal coordinate systems in the Frobenius neighborhoods of points
of $X$. By definition, $\M_{coord}(X)$ is a torsor over $X$ with
respect to the subgroup group $\Aut^0(\barA) \subset \Aut(\barA)$ of
automorphisms of the algebra $\barA$ which preserve the maximal
ideal $\m \subset \barA$. The larger group $\Aut(\barA)$ also acts
on $\M_{coord}(X)$, and the quotient is isomorphic to $X^\tw$ -- we
have a sequence of maps
$$
\begin{CD}
\M_{coord}(X) @>>> X @>{\Fr}>> X^\tw,
\end{CD}
$$
where the map on the right-hand side is the Frobenius map. Both
torsors $\M_{coord}(X)/X$ and $\M_{coord}(X)/X^\tw$ are locally
trivial in \'etale topology.

For any algebraic group $V$ equipped with an action of the group
$\Aut(\barA)$, we denote by $\Loc(V)$ the associated \'etale sheaf
on $X^\tw$. One checks easily that $\Loc(\barA) \cong \Fr_*\calo_X$,
$\Loc(k) \cong \calo_{X^\tw}$, $\Loc(\barA^*) \cong
(\Fr_*\calo_X)^*$, and $\Loc(k^*) = \calo_{X^\tw}^*$.

Just as in Proposition~\ref{tors}, one deduces from
Proposition~\ref{restr.drb} that giving a $B$-quantization of $X$ is
equivalent to giving a lifting of the torsor $\M_{coord}(X)/X^\tw$
to the group $\Aut(A^B)$ with respect to the natural group map
$$
\Aut(A^B) \to \Aut(\barA).
$$
We filter $B$ by the powers of the maximal ideal $\m \subset B$, and
we lift the torsor $\M_{coord}(X)/X^\tw$ to $\Aut(A^B)$
step-by-step, by going through the groups $\Aut(A^{B/\m^l})$. At the
first step, we have
$$
\Aut(A^{B/\m^2}) = \Aut(\barA) \rtimes \HH(\barA,\m/\m^2);
$$
therefore a lifting exists, and all liftings are classified by
elements of the cohomology group \eqref{etal}. At each subsequent
step, we apply the standard obstruction theory to \eqref{restr.ext}
and deduce that the obstruction to the lifting lies in the group
$$
H^2_{et}(X^\tw,\Loc(\HH(\barA,\m^l/\m^{l+1}))),
$$
and if this obstruction vanished, the liftings are classified by
elements of the group
$$
H^1_{et}(X^\tw,\Loc(\HH(\barA,\m^l/\m^{l+1}))).
$$
However, when $l \geq 1$, the restricted structure on
$\m^l/\m^{l+1}$ is obviously trivial. Therefore we have
$\HH(\barA,\m^l/\m^{l+1}) = (\barA/k)^N$ for some integer $N$, and
$$
\Loc(\HH(\barA,\m^l/\m^{l+1})) = (\Fr_*\calo_X/\calo_{X^\tw})^N.
$$
By assumption on $X$, this \'etale sheaf has no
cohomology. Therefore the lifting exists and is unique. This proves
the first claim. To prove the second claim, it suffices to notice
that $\Loc(A^*/k^*) = (\Fr_*\calo_X)^*/\calo_{X^\tw}^*$. Finally for
the last claim, one checks that the small Dieudonn\'e module $k
\cdot h$ is injective in the abelian category of small Dieudonn\'e
modules, so that, whatever is the restricted operation on $\m/\m^2$,
the exact sequence
$$
\begin{CD}
0 @>>> k \cdot h @>>> \m/\m^2 @>>> (\m/\m^2)/h @>>> 0
\end{CD}
$$
splits.
\endproof

When studying Frobenius-constant quantizations by
Theorem~\ref{pqua.thm}, one can further consider the Kummer spectral
sequence, and deduce the short exact sequence
\begin{equation}\label{shrt}
0 \to \Pic(X)/p\Pic(X) \to
H^1_{et}(X^\tw,(\Fr_*\calo_X)^*/\calo_{X^\tw}*) \to \Br_p(X) \to 0,
\end{equation}
where $\Pic(X^\tw) \cong \Pic(X)$ is the Picard group of $X$ and
$X^\tw$, and $\Br_p(X)$ is the $p$-torsion part of the
(cohomological) Brauer group $\Br(X^\tw) = \Br(X) =
H^2_{et}(X,\calo_X^*)$. An additional result clarifies the
appearance of the Brauer group.

\begin{prop}[{{\cite[Proposition 1.24]{bk3}}}]
In the assumption of Theorem~\ref{pqua.thm}, assume that a
Frobenius-constant quantization $\calo_h$ is classified by $a \in
H^1_{et}(X^\tw,(\Fr_*\calo_X)^*/\calo_X^{*p})$, and denote by $b \in
\Br_p(X)$ the image of the class $a$ under the canonical projection
given in \eqref{shrt}. Denote
$$
X^\tw[[h]] = \sspec(X^\tw,\calo_{X^\tw}[[h]]),\qquad\qquad
X^\tw((h)) = \sspec(X^\tw,\calo_{X^\tw}((h))),
$$
and let $\pi:X^\tw((h)) \to X^\tw$ be the natural projection.  Using
the splitting map $s:\calo_X^p \to \calo_h$, consider $\calo_h$ as a
sheaf of algebras over $X^\tw[[h]]$, and consider the localization
$\calo_h(h^{-1})$ as a sheaf of algebras over $X^\tw((h))$. Then
$\calo_h(h^{-1})$ is an Azumaya algebra over $X^\tw((h))$, and its
class in the Brauer group of $X^\tw((h))$ is equal to
$\pi^*(b)$.\endproof
\end{prop}

\section{Derived equivalence.}

Quantization theory in positive characteristic summarized in
Theorem~\ref{pqua.thm} may have some independent interest; however,
the reason for its development in \cite{bk3} was a somewhat
unexpected application purely in characteristic $0$. This is what we
are going to describe now.

\subsection{Tilting generators.}
We start with some generalities. Assume given an affine normal
algebraic variety $Y$ over a field $K$ and a smooth projective
resolution $\pi:X \to Y$. Consider the derived category $D^b_c(X)$
of bounded complexes of coherent sheaves on $X$, and assume given a
vector bundle $\E$ on $X$. Then, if we denote $R = \End(\E)$, we
have a natural functor
$$
\Rhom(\E,-):D^b_c(X) \to D^b(R\fmod),
$$
where $D^b(R\fmod)$ is the bounded derived category of finitely
generated left $R$-modules. Denote this functor by
$R^\hdot\pi_*^{\E}$. If we pass to the derived categories bounded
from above, we have an adjoint functor
$$
L^\hdot\pi^*_{\E}:D^-(R\fmod) \to D^-_c(X).
$$

\begin{defn}
\begin{enumerate}
\item A vector bundle $\E$ on $X$ is said to be {\em tilting} if
  we have $\Ext^i(\E,\E) = 0$ for all $i \geq 1$.
\item A vector bundle $\E$ is called a {\em tilting generator} for
  $X$ if in addition $\Rhom^\hdot(\E,\F)=0$ implies $\F=0$ for any
  $\F \in D^-_c(X)$.
\end{enumerate}
\end{defn}

A vector bundle $\E$ is tilting if and only if the composition
$R^\hdot\pi^{\E}_* \circ L^\hdot\pi^*_{\E}$ is the identity
endofunctor of the category $D^-(R\fmod)$. This in turn happens if
and only if the functor $L^\hdot\pi^*_{\E}:D^-(R\fmod) \to D^-_c(X)$
is fully faithful. If $\E$ is a tilting generator, then
$L^\hdot\pi^*_{\E} \circ R^\hdot\pi^{\E}_*$ is also the identity,
and both $L^\hdot\pi^*_{\E}$ and $R^\hdot\pi_*^{\E}$ are equivalence
of categories which induce equivalences between $D^b_c(X)$ and
$D^b(R\fmod)$ (for details, see \cite[Lemma 1.2]{ka1}).

Thus if $X$ admits a tilting generator $\E$, the ``geometric''
category $D^b_c(X)$ is equivalent to the purely algebraic category
$D^b(R\fmod)$. This is interesting in its own right, and also has
rather strong implications concerning the topology of $X$ which we
will describe in the next Section.

Unfortunately, tilting generators seem to be quite rare. Presently
there are only two situations where the existence of a tilting
generator is known.
\begin{enumerate}
\item $\dim X = 3$, $X$ is a crepant resolution of a quotient
  singularity $Y=V/G$, where $V$ is a $3$-dimensional vector space,
  and $G \subset SL(V)$ is a finite subgroup. This is the situation
  of a so-called {\em McKay equivalence} established in \cite{bkr}.
\item Again $\dim X = 3$, $X$ has a trivial canonical bundle, and
  $\pi:X \to Y$ is {\em small}, that is, $X$ has relative dimension
  at most $1$ over $Y$. This has been described by M. Van den Bergh
  \cite{vdb}, following the work of T. Bridgeland \cite{br}.
\end{enumerate}
It turns out that the quantization in positive characteristic allows
one to construct a tilting generator in a third rather general
situation -- namely, for a symplectic resolution $X$.

\begin{theorem}[{{\cite[Theorem 1.4]{ka1}}}]\label{equi.thm}
Let $X \to Y$ be a symplectic resolution over a field $K$ of
characteristic $0$. Then for any closed $y \in Y$, there exists an
\'etale neighborhood $Y_0 \to Y$ of the point $y \in Y$ such that
the fibered product $X_0 = Y_0 \times_Y X$ admits a tilting
generator.
\end{theorem}

Very sketchily, the reason Theorem~\ref{pqua.thm} is useful for this
result is the following. Note that if $X$ and $Y$ were defined over
a perfect field $k$ of positive characteristic $\cchar k = p > 0$,
and $H^i(X,\calo_X)=0$ for $i \geq 1$, then for any $a \in
\Pic(X)/p\Pic(X)$ Theorem~\ref{pqua.thm} gives a quantization
$\calo_a$ of $X$ associated to the image of the class $a$ in the
group $H^1_{et}(X,\calo_X^*/\calo_X^{p*})$. Moreover, by
Proposition~\ref{azu}, the sheaf $\calo_a(h^{-1})$ is a {\em split}
Azumaya algebra over $X^\tw((h))$ -- in other words,
$\calo_a(h^{-1}) = \eend(\E_a)$ for some vector bundle $\E_a$ on
$X^\tw((h))$. On the other hand, $H^i(X,\calo_a)=0$ for $i \geq 1$
by semicontinuity, and therefore $H^i(X^\tw((h)),\eend(\E_a))=0$ for
$i \geq 1$. Thus the vector bundle $\E_a$ is automatically tilting.

Elementary obstruction theory shows that tilting vector bundles are
rigid, that is, they extend uniquely to any formal deformation of
$X$. Thus by standard reduction to positive characteristic, we
obtain, in the assumption of Theorem~\ref{equi.thm}, a series of
tilting vector bundles $\E$ on $X((h)) \cong X^\tw((h))$ depending
on a prime number $p$ and a class $a \in \Pic(X)$ (since $Y$ is not
proper, we may have to replace it by a \'etale neighborhood
$Y_0$). At the cost of shrinking $Y_0$ even further, we can get of
rid of the quantization parameter $h$, and obtain a series of
tilting vector bundles on the symplectic resolution $X_0$.

Now, a careful analysis shows that for almost all values of the
parameter $a \in \Pic(X)$ the corresponding vector bundle $\E_a$
will actually be a tilting generator (in fact, it suffices to
consider values proportional to a class $[L] \in \Pic(X)$ of an
ample line bundle $L$). More precisely, there is constant $M$
independent of $p$ such that $\E_a$ is {\em not} a generator for at
most $M$ values of $a$. Thus if we take $p$ large and $a$ generic
enough, we get a tilting generator required by
Theorem~\ref{equi.thm}.

\medskip

We will now try to fill the gaps in this sketch to some extent, so
that the reader would get an idea about how the actual proof of
Theorem~\ref{equi.thm} works.

\subsection{Twistor deformations.}
The first thing to do is to collect various quantization $\E_a$ into
a single multi-parameter family. To do this, we will apply the last
claim of Theorem~\ref{pqua.thm} to a certain special one-parameter
deformation of the symplectic resolution $X$.

\medskip

Recall that if a smooth symplectic variety $Z$ over a field $k$ is
equipped with a symplectic action of the multiplicative group $\gm$,
then a map $\mu:Z \to \A^1 = \Spec k[t]$ is called a {\em moment
map} if $\Omega_Z \cntrct \xi_0 = \mu^*dt$, where $\xi_0$ is the
infinitesimal generator of the $\gm$-action. If a moment map is
given, and the quotient $Z/\gm$ exists, this quotient becomes
naturally a Poisson scheme over $\A^1$, and its fiber $X$ over the
origin $o \subset \A^1$ is again symplectic. This is known as {\em
Hamiltonian reduction}. 

It turns out that sometimes the Hamiltonian reduction procedure can
be inverted to some extent. Namely, let $X/Y$ be a symplectic
resolution over a field $k$, and let $L$ be a line bundle on
$X$. Denote $S=\Spec k[[t]]$, the formal disc over $k$, and let $o
\in S$ be the special point (given by the maximal ideal $tk[[t]]
\subset k[[t]]$).

\begin{defn}\label{tw.defn}
By a {\em twistor deformation} $Z$ associated to the pair $\langle
X,L \rangle$ we will understand a smooth symplectic deformation
$\langle \X,\LL \rangle$ of the pair $\langle X,L \rangle$ over $S$
and a symplectic form $\Omega_Z$ on the total space $Z$ of the
$\gm$-torsor associated to $\LL$ such that $\Omega_Z$ is
$\gm$-invariant, and the projection $\rho:Z \to \X \to S$ is the
moment map for the $\gm$-action on $Z$. A twistor deformation is
called {\em exact} if the symplectic form $\Omega_Z$ is exact.
\end{defn}

\begin{lemma}[{{\cite[Lemma 2.2]{ka2}}}]\label{tw.lemma}
Assume that $\cchar k =0$. Then for any line bundle $L$ on $X$,
there exists a twistor deformation $\langle \X,\LL,\Omega_Z\rangle$
associated to the pair $\langle X,L \rangle$. Moreover, $\X$ is
projective over $\Y = \Spec H^0(\X,\calo_\X)$, while $\Y$ is normal
and flat over $S$.\endproof
\end{lemma}

In terms of the period map as in Theorem~\ref{defo}, a twistor
deformation is parametrized by a straight line $[\Omega]+t[L] \in
H^2_{DR}(X)$ -- it is easy to see that this is equivalent to the
moment map condition in Definition~\ref{tw.defn}. Thus
Lemma~\ref{tw.lemma} can be deduced from
Theorem~\ref{defo}. Interestingly, twistor deformations exist in a
much wider generality -- namely, for a arbitrary Poisson scheme $X$
with $H^1(X,\calo_X)=H^2(X,\calo_X)=0$ (see \cite{ka2}).

\begin{remark}
The name {\em twistor deformation} comes from hyperk\"ahler
geometry. Namely, if the Chern class $[L] \in H^2_{DR}(X)$ can be
represented by the K\"ahler form of a hyperk\"ahler metric on $X$,
then the corresponding twistor deformation $\X/S$ can be extended
over the whole affine line $\A^1$ and in fact even over the
projective line $\Pp^1$. The total space of this extended
deformation is known as the {\em twistor space} associated to the
hyperk\"ahler metric.
\end{remark}


\begin{lemma}[{{\cite[Lemma 2.5]{ka2}}}]\label{1-1}
Assume that the line bundle $L$ on $X$ is ample, and consider the
twistor deformation $\langle \X,\LL,\Omega_Z\rangle$ associated to
the pair $\langle X,L \rangle$ by Lemma~\ref{tw.lemma}. Let $\wh{A}
= H^0(\X,\calo_{\X})$, $\Y = \Spec \wh{A}$, and let $\pi:\X \to \Y$
be the natural map. Then the map $\pi$ is one-to-one over the
complement $S \setminus \{o\}$. Moreover, if $Y$ is the spectrum of
a Henselian local $k$-algebra, so that $\wh{A}$ is a local
$k$-algebra with maximal ideal $\m \subset \wh{A}$, then there
exists a finitely generated $k$-subalgebra $\wt{A} \subset \wh{A}$
such that
\begin{enumerate}
\item the $t$-adic completion of the Henselization of the algebra
  $\wt{A}$ in $\m \cap \wt{A} \subset \wt{A}$ coincides with
  $\wh{A}$, and
\item all the data $\langle \X,\LL,\Omega_Z \rangle$ are defined
over $\wt{A}$.
\end{enumerate}
\end{lemma}

\proof[Sketch of the proof.] Since the map $\pi:\X \to \Y$ is
projective and $\Y$ is normal, for the first claim it suffices to
show that the generic fiber $\X_\eta$ over the generic point $\eta
\in S$ has no compact curves. Indeed, if $\rho:C_\eta \to \X_\eta$
is a map from a proper curve $C_\eta$, the class $[\Omega_\eta]$ of
the symplectic from $\Omega_\eta$ on $X_\eta$ satisfies
$$
\rho^*[\Omega_\eta] = \rho^*[\Omega]+t\rho^*[L],
$$
where $[\Omega]$ is the class of the form $\Omega$ on $X$, and $[L]$
is the first Chern class of the line bundle $[L]$. Since $C_\eta$ is
a curve, $\rho^*\Omega_\eta=0$, so that in particular,
$\rho^*[L]=0$. Since $L$ is ample, this implies that $\rho(C_\eta)
\subset X_\eta$ is a point, which proves the first claim. 

As a corollary, we deduce that the generic fiber $\Y_\eta \cong
X_\eta$ is smooth. Then the second claim is a particular case of
Artin's Theorem \cite[Theorem 3.9]{A}. Namely, since the formal
scheme $\Y$ contains a closed subscheme $Y \subset \Y$ if finite
type, and the complement $\Y \setminus Y \cong \Y_\eta$ is
regular, the formal scheme $\Y$ is a completion of a scheme $\wt{Y}
= \Spec \wt{A}$ of finite type.
\endproof

Assume now that we are in the situation of Theorem~\ref{equi.thm} --
we are given a symplectic resolution $X/Y$, $\pi:X \to Y$ over a
field $K$ of characteristic $0$, and a closed point $y \in Y$. By
Theorem~\ref{topo} the form $\Omega$ is exact in the formal
neighborhood of the fiber $\pi^{-1}(y) \subset X$. Changing $\Omega$
if necessary and possibly replacing $Y$ with an open neighborhood of
$y$, we can assume that $\Omega$ is exact on the whole of $X$.

Fix an ample line bundle $L$ on $X$, and consider the twistor
deformation $\langle \X,Z,\Omega_Z \rangle$ associated to $\langle
X,L \rangle$ by Lemma~\ref{tw.lemma}. One shows -- this is
\cite[Lemma 2.4]{ka1} -- that since $\Omega$ is exact on $X$, the
symplectic form $\Omega_Z$ on $Z$ is also exact (that is, the
twistor deformation is exact in the sense of
Definition~\ref{tw.defn}). Moreover, by Lemma~\ref{1-1} $\X$ and all
the other data are actually defined over a scheme $\wt{Y} = \Spec
\wt{A}$ of finite type over $K$ whose completion near $Y \subset
\wt{Y}$ is $\Y$. Therefore we can find a subring $O \subset K$ of
finite type over $\Z$ so that everything is defined and smooth over
$O$. Since $O$ is of finite type over $\Z$, the residue field $O/\m$
for any maximal ideal $\m \subset O$ is a finite field.

To sum up: starting from a symplectic resolution $X/Y$ over $K$ and
a point $y \in Y$, we can construct a symplectic resolution
$X_O/Y_O$, its exact twistor deformation $\X_O/\Y_O$, and an
$O$-valued point $y_O \in Y_O$. Localizing $O$ is necessary, we can
also assume that $\X_O$ satisfies all the topological properties of
$X$ given in Theorem~\ref{topo} -- in particular, $X_O/Y_O$ is
semismall, and $H^i(X_O,\calo_{X_O})=0$ for $i \geq 1$.

\subsection{Estimates for critical lines.}
Let us now fix $O \subset K$, the symplectic resolution $X_O/Y_O$
and its exact twistor deformation described above. For any maximal
ideal $\m \subset O$, we obtain a symplectic resolution $X_\m/Y_\m$
and its exact twistor deformation $\X_\m$ over the field $k=O/\m$ of
some positive characteristic $\cchar k = p > 0$. We also obtain the
closed point $y_\m \in Y_\m$.

Since the twistor deformation $\X_\m$ is exact, the corresponding
variety $Z_\m$ carries a restricted Poisson structure in the sense of
Definition~\ref{restr.defn}, see Proposition~\ref{restr.prop}. A
more careful analysis (\cite[Proposition 2.6]{ka1}) shows that this
restricted structure is $\gm$-invariant, so that it descends to
$\X_\m = Z_\m/\gm$, and that the deformation parameter $t$ satisfies
$t^{[p]}=t$. Setting $t^{[p]}=t$ obviously defines uniquely a
Frobenius-derivation of $k[[t]] = \calo_S$, in other words, a
restricted structure on $k[[t]]$. Thus $\X_\m/S$ is a
$k[[t]]$-quantization of the restricted Poisson variety $X_\m$.

In order to apply Theorem~\ref{pqua.thm} to $\X_\m/S$, we need to
choose a good quantization base in the sense of
Definition~\ref{qua.base.defn} which is compatible with this
restricted structure on $k[[t]]$. We let $B = k[[h,t]]$, the algebra
of power series in two variables, and define the splitting map
$s:k[[t]] \to k[[h,t]]$ by $s(t) = t(t^{p-1}-h^{p-1})$. Then
$\langle B,s \rangle$ obviously satisfies the conditions of
Definition~\ref{qua.base.defn}, and
$$
t^{[p]} = \frac{1}{h^{p-1}}(t^p-s(t)) = t,
$$
as required. Applying the last claim of Theorem~\ref{pqua.thm}, we
obtain a $B$-quantization $\calo_h$ of the Poisson scheme $X_\m$
extending the given $B/h$-quantization $\X_\m$ -- in other words, a
Frobenius-constant quantization $\calo_h$ of the restricted Poisson
scheme $\X_\m$ which is a sheaf of algebras over $B$.

Geometrically, the map $s:S_h = \Spec B \to S$ is a flat map with a
singular fiber over the origin $o \in S$. This fiber $S_o=s^{-1}(o)
\subset S_h$ is the union of the lines $S_a \subset S_h$ defined by
equations $h-at$ for all $a \in \Z/p\Z$. We let $S_t = \Spec
k[t]/t^p = S_o \times_{k[[h]]} k$. By definition, $\calo_h$ is a
sheaf of algebras on
$$
\X_h = \X^\tw_\m \whtimes_{S^\tw} S_h,
$$
where $\whtimes$ means the tensor product completed with respect to
the natural adic topology.  For every line $S_a \subset S_h$, the
subscheme $X_a = \X^\tw \whtimes_{S^\tw} S_a \subset \X_h$ is
canonically identified with $X_h = X^\tw_\m[[h]]$.  If we restrict
$\calo_h$ to $X_a \subset \X_h$, we obtain a Frobenius-constant
quantization $\calo_a$ of the scheme $X_\m$ associated by
Theorem~\ref{pqua.thm} to the parameter $a[L] \in
\Pic(X_\m)/p\Pic(X_\m) \subset
H^1_{et}(X_\m,\calo_{X_\m}^*/\calo_{X_\m}^{*p})$.

By Proposition~\ref{azu}, the algebra $\calo_a(h^{-1})$ is a matrix
algebra for any $a \in \Z/p\Z$, that is,
$\calo_a(h^{-1})=\eend(\E_a)$ for some vector bundle $\E_a$ on
$\overline{X}_a = X_a \setminus X_\m^\tw \cong X_\m((h))$. If we
consider the union
$$
X_o = \bigcup_a X_a = X^\tw_\m \whtimes S_o \subset \X_h,
$$
then the complement $\overline{X}_o = X_o \setminus (X^\tw_\m \times
S_t)$ is the disjoint union of the complements
$\overline{X}_a$, and $\calo_h$ restricts to a matrix algebra over
the whole $\overline{X}_o$.

\begin{defn}\label{reg.line}
A line $S_a \subset S_h$ is said to be {\em regular} for $X_\m$ if
the vector bundle $\E_a$ is a tilting generator on the completion of
$X^\tw_\m((h))$ near the preimage $\pi^{-1}(y_\m((h)))$. A line which is
not regular is called {\em critical}.
\end{defn}

The main technical result needed for Theorem~\ref{equi.thm} is the
following.

\begin{prop}\label{est}
There exists a constant $C$, independent of the choice of the ideal
$\m \subset O$, such that there are at most $C$ lines $S_a \subset
S_h$ which are critical for $X_\m$.
\end{prop}

Thus if we take the ideal $\m \subset O$ so that $p=\cchar O/\m$ is
high enough, there exists at least one regular line $S_a$, $a \in
\Z/p\Z$.

\medskip

The proof of Proposition~\ref{est} takes up Section 3 and most of
Section 4 of \cite{ka1}; we refer a reader interested in technical
details to that paper. Here we only list the main steps of the
proof.

\step{} We need a way to measure whether the direct image functor
$R^\hdot\pi_*^\E$ is an equivalence. We work in a general situation:
consider a scheme $X$ equipped with a coherent sheaf of algebras
$\cala$ such that $H^i(X,\cala)=0$ for $i \geq 1$, and consider the
global sections functor $R^\hdot\pi_*^\cala:D^b_c(X,\cala) \to
D^b(A\fmod)$ from the derived category of coherent sheaves of
$\cala$-modules to the category of finitely generated modules over
$A = H^0(X,\cala)$. The functor has an adjoint
$L^\hdot\pi^*_\cala:D^-(A\fmod) \to D^-_c(X,\cala)$, and since $A =
R^\hdot\pi^\cala_*(\cala)$, we have $R^\hdot\pi^\cala_* \circ
L^\hdot\pi^*_\cala \cong \Id$. We show (see \cite[Subsection
3.1]{ka1}) that the composition in the other direction, -- that is,
the functor
\begin{equation}\label{comp.eq}
L^\hdot\pi^*_\cala \circ R^\hdot\pi_*^\cala:D^-_c(X,\cala) \to
D^-_c(X,\cala)
\end{equation}
is a kernel functor associated to a kernel $\M^\hdot(X,\cala) \in
D^-_c(X \times X,\cala^{opp} \boxtimes \cala)$. The kernel
$\M^\hdot(X,\cala)$ maps naturally to the diagonal sheaf
$\cala_\Delta$ on $X \times X$, and the cone $\K^\hdot(X,\cala)$ of
this map is trivial if and only if the functor \eqref{comp.eq} is an
equivalence. The kernel $\K^\hdot(X,\cala)$ essentially depends only
on the topology of $X$, not on its scheme structure; in particular,
for any finite radical map $\rho:X \to X'$ we have
\begin{equation}\label{radi}
\K^\hdot(X',\rho_*\cala) \cong \rho_*(X,\cala).
\end{equation}
If the scheme $X$ is flat over a scheme $S$, then a
similar construction exists relatively over $X$, and the kernel
$\K^\hdot(X,\cala)$ is actually supported on $X \times_S X \subset X
\times X$. Moreover, this construction is compatible with the base
change: for any map $\rho:S' \to S$, we have
$\rho^*\K^\hdot(X,\cala) \cong \K^\hdot(X \times_S S',
\rho^*\cala)$. Finally, we note that if $D^b_c(X,\cala)$ has finite
homological dimension, say $d$, then it suffices to check that
$\K^i(X,\cala)$ is trivial for $0 \geq i \geq -2d$ -- if this is
satisfied, then $\K^i(X,\cala)=0$ for all $i$.

\step{} Applying Step 1 to Proposition~\ref{est}, we see that a line
$S_a \subset S_h$ is regular if and only if the kernel
$$
\overline{\K}^\hdot_a = \K^\hdot(\overline{X}_a,\calo_a) =
\K^\hdot(X^\tw_\m((h)),\eend(\E_a))
$$
is trivial on the fiber $F_y=\pi^{-1}((y_\m \times y_\m)((h)))
\subset (X^\tw_\m \times X^\tw_\m)((h))$. Moreover, since
$D^b_c(X^\tw((h)))$ has homological dimension $\dim X$, we can claim
that if the line $S_a$ is critical, then there exists a non-positive
integer $i \geq -2\dim X$ and a point $x \in F_y$ such that
$\K^j_a=0$ on $F_y$ for all $j > i$, while the fiber
$\left(\K^i_a\right)_x$ is non-trivial. However, $\K^i_a$ is by
definition a sheaf of modules over $\eend(\E_a)^{opp} \boxtimes
\eend(\E_a)$, and this is a matrix algebra of rank $p^{2\dim
X}$. Therefore we can claim that
\begin{equation}\label{est.1}
\dim_{k((h))}\left(\overline{\K}^i_a\right)_x \geq p^{\dim X}.
\end{equation}

\step{} Choose a projective embedding $X^\tw_\m \times X^\tw_\m \to
\PP_Y = \Pp^N \times Y^\tw_\m \times Y^\tw_\m$, and denote by $\PP_y
\subset \PP_Y$ the fiber of $\PP_Y/Y$ over the point $y_\m \times
y_\m \in Y^\tw_\m \times Y^\tw_\m$, with the embedding map
$i_y:\PP_y \to \PP_Y$. Let $\PP_h = \PP_y((h)) \subset
\PP_Y((h))$. Applying pushforward, we can treat the complex
$\overline{\K}^\hdot_a$ as a complex of sheaves on $\PP_Y$; by
restriction we obtain a complex of sheaves
$\overline{\K}^\hdot_{y,a} = L^\hdot i_y^*\overline{\K}^\hdot_a$ on
$\PP_h$, and \eqref{est.1} implies
$$
\dim_{k((h))}\Hh^i\left(\PP_h,\overline{\K}^\hdot_{y,a}
\overset{L}{\otimes} k_x\right) =
\dim_{k((h))}\left(\overline{\K}^i_{y,a}\right)_x \geq p^{\dim X}.
$$ 
But by definition, $\PP_h$ is the projective space $\Pp^N_{k((h))}$
of dimension $N$ over $k((h))$. Therefore the skyscraper sheaf $k_x$
at the point $x \in \PP_h$ has a Koszul resolution by sheaves of the
form $\calo(-n)^{\oplus{\binom{N}{n}}}$, $0 \geq n \geq N$, and we
conclude that
$$
p^{\dim X} \leq \dim_{k((h))}\Hh^i\left(\PP_h,\overline{\K}^\hdot_{y,a}
\overset{L}{\otimes} k_x\right) \leq \sum_{0 \leq n \leq N}
\dim_{k((h))}\binom{N}{n}\Hh^{i-n}(\PP_h,\overline{\K}^\hdot_{y,a}(-n)).
$$
The right-hand side does not depend on the point $x \in \PP_h$.
Therefore we can sum up these estimates over all critical $a \in
\Z/p\Z$, and conclude that to bound from above the number of
critical lines, it suffices to prove that for any $n$, $0 \leq n
\leq N$ and $i$, $0 \geq i \geq -2\dim X$, we have
\begin{equation}\label{est.2}
\sum_{a \in \Z/p\Z}\dim_{k((h))}
\Hh^{i-n}(\PP_h,\overline{\K}^\hdot_{y,a}(-n)) \leq M p^{\dim X},
\end{equation}
where $M$ is a positive integer which does not depend on the maximal
ideal $\m \in O$.

\step{} Using the disjoint union decomposition $\overline{X}_o =
\coprod_a \overline{X}_a$, we can rewrite \eqref{est.2} as
$$
\dim_{k((h))}\Hh^{i-n}(\PP_o,(L^\hdot
i_y^*\overline{\K}^\hdot_o)(-n)) \leq M p^{\dim X},
$$
where $\PP_o = \PP_y \whtimes S_o \subset \PP_Y \whtimes S_o$, and
$\overline{\K}^\hdot_o =
\K^\hdot(\overline{X}_o,\calo_h)$. Moreover, since the scheme $S_o$
is flat over $k[[h]]$ and formation of the kernel $\K^\hdot(-,-)$
commutes with base change, we can apply semicontinuity theorem to
$X_o/k[[h]]$; we conclude that to prove \eqref{est.2}, it suffices
to prove that
\begin{equation}\label{est.3}
\dim_k\Hh^{i-n}(\PP_y \times S_t,(L^\hdot
i_y^*\K^\hdot_t)(-n)) \leq M p^{\dim X},
\end{equation}
where
$$
\K^\hdot_t = \K^\hdot(X_o/h,\calo_o/h) =
\K^\hdot(X_\m^\tw \times S_t,\calo_{\X}/t^p).
$$

\step{} It remains to prove \eqref{est.3}. It explicitly does not
depend on the quantization $\calo_h$, but still depends on the
twistor deformation $\X_\m$. We first note that the kernel
$\K^\hdot_t$ can be rewritten as
$$
\K^\hdot_t = \K^\hdot(\X_\m^{[p]},\Fr^S_*\calo_{\X_\m})/t^p,
$$
where $\X_\m^{[p]} = \X_\m^\tw \times_{S^\tw} S \subset
\X_\m$, and $\Fr^S$ is the embedding map (simultaneously, the relative
Frobenius map for $\X/S$). But since the map $\Fr^S$ is finite and
radical, we may apply \eqref{radi} and write
$$
\K^\hdot_t = \Fr^S_*\K^\hdot/t^p,
$$
where $\K^\hdot = \K^\hdot(X_\m,\calo_{X_\m})$. Since the map
$\pi:\X_\m \to \Y_\m$ is one-to-one over the generic point of $S$,
the kernel $\K^\hdot$ is annihilated by $t^{M_1}$ for some integer
$M_1$. Moreover, the map $\pi$ is generically one-to-one already
over our original field $K$ of characteristic $0$, so that $M_1$
does {\em not} depend on the choice of $\m \subset O$. We deduce
that
$$
\dim_k\Hh^{i-n}(\PP_y \times S_t,L^\hdot
i_y^*\K^\hdot_t(-n)) \leq M_1\dim_k\Hh^{i-n}(\PP_y,L^\hdot
i_y^*(\Fr_*\K^\hdot)(-n)),
$$
where $M_1$ does not depend on $\m$. Thus to prove \eqref{est.3}, it
suffice to prove that
\begin{equation}\label{est.4}
\dim_k\Hh^{i-n}(\PP_y,L^\hdot i_y^*(\Fr_*\K^\hdot)(-n)) \leq M
p^{\dim X}
\end{equation}
for some universal constant $M$.

\step{} Finally, \eqref{est.4} only depends on $X_\m/Y_\m$, not on
the quantization $\calo_h$, nor on the twistor deformation
$\X$. Thus the symplectic form on $X$ is no longer used. We can
compactify $X/Y$ to a birational map $\wt{X} \to \wt{Y}$ of proper
schemes over $K$, obtain $\wt{X}_O$ and $\wt{X}_\m$ (possibly
changing $O \subset K$), and extend $\K^\hdot_O =
\K^\hdot(X_O,\calo_{X_O})$ to some complex $\wt{\K}_O^\hdot$ of
coherent sheaves on $\wt{X}_O \times \wt{X}_O$. We denote its
restriction to $X \times X \subset X_O \times_O X_O$ by
$\wt{\K}^\hdot$, and we denote its restriction to $X_\m \times X_\m
\subset X_O \times_O X_O$ by $\wt{\K}^\hdot_\m$. Choose a projective
embedding $\wt{X}_O \times_O \wt{X}_O \to \PP^{N_1}_O$. Localizing
$O$ if necessary, we may assume that the Hilbert polynomials of the
sheaves $\wt{\K}^\hdot_\m$ do not depend on $\m \subset O$ and
coincide with the Hilbert polynomials of the sheaves
$\wt{\K}^\hdot$.

We can now remove the restriction functor $L^\hdot i_y^*$ from
\eqref{est.4}. To do this, we apply again the Koszul resolution
argument of Step 3, and conclude that it suffices to find an
estimate of the form
\begin{equation}\label{est.5}
\dim_k\Hh^{i-n}(\wt{X}^\tw_\m \times
\wt{X}^\tw_\m,(\Fr_*\wt{\K}^\hdot)(-n)) \leq M p^{\dim X}
\end{equation}
But the left-hand side is equal to
$$
\dim_k\Hh^{i-n}(\wt{X}_\m \times \wt{X}_\m,\wt{\K}^\hdot(-pn)),
$$
so that \eqref{est.5} is a statement about the Hilbert polynomials
of the sheaves $\wt{\K}^\hdot$ -- which are by assumption the same
as the Hilbert polynomials of the sheaves $\wt{\K}^\hdot$. We have
to show that the degrees of these polynomials are not greater than
$\dim X$. Indeed, $\K^\hdot$ is supported on $X \times_Y X \subset X
\times X$, and choosing the extension $\wt{\K}^\hdot$ in an
appropriate way, we can also insure that $\dim\Supp\wt{\K}^\hdot =
\dim\Supp\K^\hdot \leq \dim X \times_Y X$. Since $X/Y$ is a
symplectic resolution, it is semismall by Theorem~\ref{topo}, in
other words, $\dim X \times_Y X = \dim X$. Thus
$\dim\Supp\wt{\K}^\hdot \leq \dim X$, which yields \eqref{est.5}.

\subsection{Artin approximation.} Using Proposition~\ref{est}, the
proof of Theorem~\ref{equi.thm} proceeds as follows. Assume given a
symplectic resolution $X/Y$ over a field $K$ of characteristic $0$
and a point $y \in Y$. As in the last Subsection, we choose a
subalgebra $O \subset K$ of finite type over $\Z$, schemes $X_O/Y_O$
flat, smooth and of finite type over $\Spec O$, and an $O$-valued
point $y_O:\Spec O \to Y_O$ such that $X_O \otimes_O K \cong X$,
$Y_O \otimes_O K \cong Y$, $y_O \otimes_O K = y$, and $X_O$ is
projective over $Y_O$. For any maximal ideal $\m \subset O$, we
obtain by reduction schemes $X_\m/Y_\m$ and a point $y_\m \in
Y_\m$. By Proposition~\ref{est}, we can choose $\m \subset O$ in
such a way that $X_\m$ admits a regular line in the sense of
Definition~\ref{reg.line}. Explicitly, consider the point $y_\m((h))
\in Y_\m((h))$, let $\wh{Y_\m}$ be the completion of $Y_\m((h))$
near $y_\m((h))$, and let $\wh{X_\m} = \wh{Y_\m} \times_{Y_\m}
X_\m$; then by Proposition~\ref{est} we have a vector bundle $\E_\m
= \E_a$ on $\wh{X_\m}$ which is a tilting generator.

Let $\wh{Y_O}$ be the completion of $Y_O \otimes_O O((h))$ near
$y_\m((h))$, and let $\wh{X_O} = \wh{Y_O} \times_{Y_O} X_O$. Then
$\wh{X_O}$ is flat and smooth over the completion $\wh{O}$ of the
algebra $O((h))$ with respect to the maximal ideal $\m((h))$, and
the special fiber of $\wh{X_O}/\wh{X_O}$ is identified with
$\wh{X_\m}$. Since the vector bundle $\E\m$ is tilting, it extends
uniquely to $\wh{X_O}$ considered as a a formal scheme -- indeed, by
standard deformation theory obstructions to this at each level of
the adic filtration lie in $\Ext^2(\E_\m,\E_\m)$, the choices of
extensions are parametrized by elements of $\Ext^1(\E_\m,\E_\m)$,
and both groups are trivial. By \cite[Th\'eor\`eme 5.4.5]{EGA}, the
vector bundle $\E_\m$ therefore extends to a vector bundle $\E_O$
over the actual scheme $\wh{X_O}/\wh{Y_O}$. By Nakayama Lemma the
vector bundle $\E_O$ is also tilting, and the corresponding kernel
$\K^\hdot(\wh{X_O},\eend(\E_O))$ vanishes, so that it is a tilting
generator.

By Artin Approximation Theorem \cite[Theorem 1.10]{A}, there exists
a subalgebra $O' \subset \wh{O}$ of finite type over $O$, schemes
$X_{O'}/Y_{O'}$ of finite type over $O'$, and an $O'$-valued point
$y_{O'}:\Spec O' \to Y_{O'}$ such that $X_{O'} \otimes_{O'} \wh{O}
\cong \wh{X_O}$, $Y_{O'} \otimes_{O'} \wh{O} \cong \wh{Y_O}$,
$y_{O'} \otimes_{O'} \wh{O} \cong \wh{y_O}$, $\X_{O'}$ is projective
over $Y_{O'}$, $\E_O$ by approximated to a high order by $\E_{O'}
\otimes_{O'} \wh{O}$ for a vector bundle $\E_{O'}$ on $X_{O'}$, and
on the other hand, the natural maps $X_{O'} \to X_O \times_O O'$,
$Y_{O'} \to Y_O \times_O O'$ are \'etale, and the second map sends
$y_{O'}$ to $y_O \times_O O'$. Again by Nakayama Lemma, we note that
if the order of approximation is high enough, then shrinking
$Y_{O'}$ if necessary, we can guarantee that the vector bundle
$\E_{O'}$ is a tilting generator for $X_{O'}$.

It remains to take a generic point $o \in \Spec O'$ whose residue
field $K'$ is a finite extension of our original field $K$, and
notice that, possibly after shrinking $Y_{O'}$ even further, $Y_0 =
Y_{O'} \otimes_{O'} K'$ is an \'etale neighborhood of the point $y
\in Y$, and $\E = \E_{O'} \otimes_{O'} K'$ is a tilting generator
for $X_0 = X_{O'} \otimes_{O'} K' = Y_0 \times_Y X$.

\section{Geometric corollaries.}

\subsection{Additional results on derived equivalences.} One
unsatisfactory thing about Theorem~\ref{equi.thm} is the need to fix
a point $y \in Y$ and pass to an \'etale neighborhood $Y_0$. Of
course, one can cover the whole of $Y$ with such \'etale
neighborhoods, but at present, we do not know whether the tilting
generators provided by Theorem~\ref{equi.thm} patch together. There
very well might be an obstruction to this lying in the Brauer group
$\Br(X)$. In practice, this problem is alleviated by the following
additional result.

\begin{theorem}[{{\cite[Theorem 1.8]{ka1}}}]\label{gm.act}
In the assumptions of Theorem~\ref{equi.thm}, assume in addition
that $Y$ admits a positive-weight $\gm$-action. Then this action
lifts canonically to a $\gm$-action on $X$, and the tilting
generator $\E$ provided by Theorem~\ref{equi.thm} extends to a
$\gm$-equivariant tilting generator on the whole of $X$.\endproof
\end{theorem}

In fact, Theorem~\ref{gm.act} is valid for any tilting generator of
the type provided by Theorem~\ref{equi.thm}; where it came from is
not relevant to the proof. In light of Conjecture~\ref{gm.conj},
Theorem~\ref{gm.act} is potentially very useful.

By general nonsense, the presence of a tilting generator yields
strong restrictions on the topology of a resolution $X$, further
extending those given in Theorem~\ref{topo}. Namely, we have the
following.

\begin{theorem}[{{\cite[Theorem 1.9]{ka1}}}]\label{diag}
Assume that a smooth manifold $X$ is projective over an affine local
Henselian scheme $Y/k$ and admits a tilting generator $\E$. Then the
structure sheaf $\calo_\Delta$ of the diagonal $\Delta \subset X
\times X$ admits a finite resolution by vector bundles of the form
$\E_i \boxtimes \F_i$, where $\E_i$, $\F_i$ are some vector bundles
on $X$.\endproof
\end{theorem}

\begin{corr}[{{\cite[Corollary 1.10]{ka1}}}]\label{pure}
Assume that a smooth manifold $X$ is projective over an affine
scheme $Y$, and let $E \subset X$ be the fiber over a closed point
$y \in Y$. Assume that $Y$ admits a positive-weight $\gm$-action
that fixes $y \in Y$, and assume that $X$ admits a tilting generator
$\E$. Then the cohomology groups $H^\hdot(E)$ of the scheme $E$ are
generated by classes of algebraic cycles.\endproof
\end{corr}

In this Corollary we are deliberately vague as to what particular
cohomology groups $H^\hdot(E)$ one may take. In fact, every
cohomology theory with the standard weight formalism will suffice;
in particular, the statement is true for $l$-adic cohomology and for
analytic cohomology when the base field $K$ is $\C$. The proof is
rather standard: one considers the identity endomorphism of the
cohomology $H^\hdot(E)$ and, using Theorem~\ref{diag}, decomposes it
as
\begin{equation}\label{diag.eq}
\Id(a) = \sum \eta_i(a)[a_i],
\end{equation}
where $\eta_i$ are certain linear forms on $H^\hdot(E)$, and $[a_i]$
are classes of algebraic cycles. However, there is a complication --
since the scheme $X$ is not compact, the natural map from the
cohomology $H^\hdot_c(X)$ with compact support to the usual
cohomology $H^\hdot(X)$ is not at all an isomorphism, and the usual
way to deduce \eqref{diag.eq} does not work. To overcome this
difficulty, we have to require an existence of a $\gm$-action and
work with $\gm$-equivariant cohomology. This seems much too strong;
however, at present, we do not know whether Corollary~\ref{pure}
holds without the $\gm$-action assumption.

Theorem~\ref{diag} itself is a direct corollary of the equivalence
\begin{equation}\label{equiv.eq}
D^b_c(X) \cong D^b(R\fmod),
\end{equation}
where $R = \End(\E)$. The non-commutative algebra $R$ is finite over
its Henselian center, so that it has a finite number of
indecomposable projective modules $P_i$. The equivalence
\eqref{equiv.eq} sends $R$ itself to $\E$; every projective module
$P_i$, being a direct summand of $R^N$ for some $N$, goes to a
vector bundle $\F_i$ on $X$, and it is these vector bundles that
appear in the resolution of the diagonal.

The algebra $R$ also has a finite number of irreducible modules;
those go to some complexes of coherent sheaves on $X$ supported near
the exceptional fiber $E \subset X$. In fact, one can use the
equivalence \eqref{equiv.eq} to translate the standard $t$-structure
on $D^b(R\fmod)$ to a rather unusual $t$-structure on $D^b_c(X)$ --
it is this ``perverse'' $t$-structure on the category of coherent
sheaves on $X$ that was discovered in \cite{br} in $\dim 3$. The
perverse $t$-structure is Artinian and Noetherian. Its irreducible
objects provide a canonical basis in the $K$-group $K'_0(X)$. It
would be very interesting to compute this basis in various
particular cases, such as the quiver variety case
(Example~\ref{qui.exa}).

We note that in our construction of the tilting generator $\E$,
there are three choices: we have to choose a maximal ideal $\m
\subset O$ with residue field $k=O/\m$ of positive characteristic
$\cchar k = p$, an ample line bundle $L$ on $X$, and a regular value
$a \in \Z/p\Z$. Since by construction, $\E$ is a vector bundle of
rank $p^{\dim X}$, it obviously depends at least on the residual
characteristic $p$. However, we venture the following.

\begin{conj}\label{uni}
In the assumption of Theorem~\ref{diag}, the perverse $t$-structure
induced on the derived category $D^b_c(X)$ is the same, up to a
twist by an autoequivalence of $D^b_c(X)$, for almost all maximal
ideals $\m \subset O$, ample line bundles $L$ on $X$, and regular
values $a \in \Z/p\Z$.
\end{conj}

``Almost all'' here means, hopefully, ``all but a finite
number''. Unfortunately, our methods do not yield an easy way to
compare the results for different values $a \in \Z/p\Z$, and
comparison between different maximal ideals $\m \subset O$ seems
to be completely out of reach.

In the simplest possible example $X = T^*\Pp^1$, the cotangent
bundle to $\Pp^1$, one can follow through the proof of
Theorem~\ref{equi.thm} in an effective way, with the following
end result:
$$
\E_a \cong \calo_X^{\oplus a} \oplus \calo_X(1)^{\oplus (p-a)},
$$
where $\calo_X(1)$ is the pullback of the standard line bundle
$\calo(1)$ on $\Pp^1$ with respect to the projection $X = T^*\Pp^1
\to \Pp^1$. Thus every $a \neq 0$ is regular, and the tilting
generator $\E_a$ is the sum of two vector bundles $\calo_X$,
$\calo_X(1)$ with different multiplicities depending on $\m \subset
O$ and $a \in \Z/p\Z$. It is easy to see that Conjecture~\ref{uni}
is true in this case, with all the tilting generators giving the
same $t$-structure as $\calo_X \oplus \calo_X(1)$ (which is also a
tilting generator, and in a sense, the smallest possible one). We
expect that in general, the situation is the same: there is a finite
number of indecomposable vector bundles $\E_i$ which generate the
$t$-structure, and all the tilting generators $\E_a$ are obtained by
summing up the bundles $\E_i$ with multiplicities depending on $\m
\subset O$, $L$, and $a \in \Z/p\Z$.

\begin{remark}
There is in fact one more choice in the proof of
Theorem~\ref{equi.thm} which we tacitly ignore in the above
discussion: when we represent the matrix algebra $\calo_a(h^{-1})$
as an endomorphism algebra $\eend(\E_a)$, the vector bundle $\E_a$
is only defined by up to a twist by a line bundle. The ``twist by an
autoequivalence'' clause in Conjecture~\ref{uni} is needed to take
care of this. To be on the safe side, we do not require this
autoequivalence to be a twist by a line bundle. In general, it would
be very interesting to study the group of all autoequivalences of
the triangulated category $D^b_c(X)$ and its action on various
perverse $t$-structures; however, at present there seems to be no
way to do this, at least in the interesting case $\dim X > 2$.
\end{remark}

One additional observation is the following.

\begin{prop}\label{defo.tilt}
In the assumptions of Theorem~\ref{equi.thm}, every tilting
generator $\E$ on $X$ extends uniquely to a tilting generator
$\wt{\E}$ on the universal deformation $\X$ provided by
Theorem~\ref{defo}.
\end{prop}

\proof{} Standard deformations theory: since $\E$ is tilting,
$\Ext^i(\E,\E) = 0$ for $i=1,2$; thus there are no obstructions to
deforming it together with $X$, and no choices are involved in such
a deformation. By Nakayama Lemma, the deformed vector bundle
$\wt{\E}$ is also a tilting generator.  \endproof

Thus Theorem~\ref{gm.act}, Theorem~\ref{diag}, Corollary~\ref{pure},
and all the above discussion apply just as well to the scheme $\X$.

\begin{remark}
In the case $\dim X = 3$, $K_X$ trivial -- that is, in the case
studied in \cite{br} and \cite{vdb} -- it is known that $Y$, being a
terminal singularity, must be a so-called cDV point, and the whole
$X/Y$ is a one-parameter deformation of a partial resolution
$X_0/Y_0$ of a Du Val quotient singularity $Y_0 = \A^2/G$, $G
\subset SL(2,K)$. However, $X_0$ is usually singular -- it is only
the total space $X$ of the deformation that is smooth. Thus $X/Y$ is
not really of the form $\X/\Y$ for some $2$-dimensional symplectic
resolution, and our methods do not apply. It would be very
interesting to try to generalize our approach to this situation and
compare it with \cite{vdb}.
\end{remark}

Finally, there is a result which compares the derived categories
$D^b_c(X)$ for different crepant resolutions of the same symplectic
singularity $Y$. This is a generalization of the particular case of
\cite[Section 3, Conjecture]{BO1} proved by Y. Kawamata in
\cite{Ka}: in Kawamata's language, ``$K$-equivalence implies
$D$-equivalence''. We would also like to mention that a particular
case of this result was proved by Y. Namikawa in \cite{N3}.

\begin{theorem}[{{\cite[Theorem 1.6]{ka1}}}]
In the assumptions of Theorem~\ref{equi.thm}, assume given a
different resolution $X'$, $\pi':X \to Y$ of the variety $Y$ with
trivial canonical bundle $K_{X'}$. Then for every closed point $y
\in Y$, there exists an \'etale neighborhood $Y_y \to Y$ such that
the derived categories $D^b_c(X \times_Y Y_y)$ and $D^b_x(X'
\times_Y Y_y)$ are equivalent.
\end{theorem}

\proof[Sketch of the proof.] One checks easily that since $K_{X'}$
is trivial, the resolution $X'$ must also be symplectic. Since a
symplectic resolution in $\dim 2$ admits a {\em unique} symplectic
resolution, the rational map $X \ratto X'$ induces an isomorphism
$X_0 = \pi^{-1}(Y_0) \cong X_0' = (\pi')^{-1}(Y_0)$, where the open
subset $Y_0 \subset Y$ is the union of the strata of dimensions
$\dim Y$ and $\dim Y - 2$ with respect to the stratification of
Theorem~\ref{sympl.singu.th}.

Going through the proof of Theorem~\ref{equi.thm}, we choose an
ample line bundle $L$ on $X$ and obtain a tilting generator $\E_a$
for $X \times_Y Y_y$; repeating the same argument for $X'$ equipped
with the strict transform $L'$ of the line bundle $L$, and possibly
changing $Y_y$, we obtain a tilting vector bundle $\E'_a$ on $X'
\times_Y Y_y$.

We can {\em not} claim that $\E'_a$ is a tilting generator: indeed,
unless $X \cong X'$, the line bundle $L'$ is not ample on
$X'$. However, since $X'/Y$ is semismall, the complements $X
\setminus X_0 \subset X$, $X' \setminus X'_0 \subset X'$ are of
codimension at least $2$. Moreover, $H^i(X_0,\calo_{X_0}) = 0$ for
$i=1,2$, and, analyzing the proof of Theorem~\ref{pqua.thm}, we
conclude that the quantizations used in the construction of the
tilting bundles $\E_1$, $\E'_a$ agree on $X_0$. Therefore $\E_a
\cong \E_a'$ on $X_0 \times_Y Y_y$. Again, since the complement to
$X_0$ is of high codimension both in $X$ and $X'$, we conclude that
the algebra $R = \End(\E_a)$ is isomorphic to $R' = \End(\E'_a)$.

In particular, the algebra $R'$ has finite homological dimension, so
that the natural functor $D^-(R'\fmod) \to D^-_c(X')$ induces a
functor
$$
D^b(R'\fmod) \to D^b_c(X').
$$
Since $\E_a'$ is tilting, this functor is a fully faithful
embedding with admissible image in the sense of \cite{BO2}. To
finish the proof, it suffices to use the following standard trick.

\begin{lemma}\label{cy}
Assume given an irreducible smooth variety $X$ with trivial
canonical bundle $K_X$ equipped with a birational projective map
$\pi:X \to Y$ to an affine variety $Y$. Then any non-trivial
admissible full triangulated subcategory in $D^b(X)$ coincides with
the whole $D^b(X)$.
\end{lemma}

For the proof we refer the reader, for instance, to \cite[Section
2]{bk2}.
\endproof

We note that in general, Lemma~\ref{cy} gives a quick and easy way
to prove that a tilting vector bundle $\E$ is a generator, avoiding
all the difficult estimates of Proposition~\ref{est}. However, in
order to apply it, one need to know that the algebra $R = \End(\E)$
has finite homological dimension. It seems that in general, there is
no way to prove it short of proving that $\E$ is a generator. 

One notable exception to this is the quotient singularity case $Y =
V/G$ considered in \cite{bk2} (this is our
Example~\ref{quo.exa}). In this case, using a more detailed analysis
of quantizations, one shows that there exists a tilting vector
bundle $\E$ on $X$ such that $\End(\E) \cong S^\hdot(V^*) \# G$, the
smash-product of the algebra of polynomial functions on $V$ and the
group algebra of the group $G$. This algebra obviously has finite
homological dimension; therefore $\E$ is a generator. No version of
Proposition~\ref{est} is required, and the description of $D^b_c(X)$
is more explicit than in the general case.

\begin{remark}
One thing that was not done in \cite{bk2} is the analysis of the
deformed tilting generator $\wt{\E}$ provided by
Proposition~\ref{defo.tilt}. The endomorphism algebra $\wt{R} =
\End(\wt{\E})$ is a flat deformation of the endomorphism algebra
$\End(\E) \cong S^\hdot(V^*) \# G$. One expects that $\wt{R}$
coincides with the so-called {\em symplectic reflection algebra}
introduced in \cite{eg}, but this has never been verified expect in
some special cases, see \cite{go}.
\end{remark}

\subsection{Comparison with quantum groups.}\label{quantu}
To finish the paper, we would like to return to the starting point
mentioned in the Introduction and give some speculations on the
connections of the present work with Geometric Representation Theory.

The motivation for the research carried out in \cite{bk2},
\cite{bk3}, \cite{ka1} was the paper \cite{BMR} and its sequel
\cite{BMR1}, where the authors study the case $X=T^*M$, $M=G/P$, a
partial flag variety associated to a semisimple algebraic group $G$
and a parabolic subgroup $P \subset G$. In that case, a particular
series of quantizations of the cotangent bundle $X=T^*M$ is given
from the very beginning -- one can consider the algebras $\D_{M,L}$
of differential operators on $M$ twisted by a line bundle $L$. The
classic result of A. Beilinson and J. Bernstein \cite{BB} claims
that in characteristic $0$, the partial flag variety $M=G/P$ is {\em
$D$-affine} for generic $L$ -- that is, the category of sheaves of
$\D_{M,L}$-modules is equivalent to the category of modules over the
algebra $H^0(M,\D_{M,L})$ of global sections of the sheaf
$\D_{M,L}$. In positive characteristic, the statement is no longer
true; however, and it has been proved in \cite{BMR}, \cite{BMR1},
the equivalence does survive on the level of derived categories:
the natural global sections functor induces an equivalence between
the derived categories $D^b_c(M,\D_{M,L})$ and
$D^b(H^0(M,\D_{M,L})\fmod)$. Moreover, the algebra $\D_{M,L}$
acquires a large center, so that sheaves of $\D_{M,L}$-modules on
$M$ can be localized to sheaves on $X^\tw$. As in our
Theorem~\ref{equi.thm}, the equivalence can then be lifted back to
characteristic $0$; the resulting algebra of global sections is
closely related to the so-called quantum enveloping algebra at a
$p$-th root of unity (see \cite{Kobi}).

In general, one can use the sheaves $\D_{M,L}$ for the cotangent
bundle $T^*M$ of any algebraic variety $M$, but this is not expected
to be very useful -- indeed, the derived $D$-affine property in
positive characteristic would in particular imply that $X = T^*M$
satisfies the assumptions of Conjecture~\ref{dema}, so that we are
automatically in the situation of \cite{BMR1}. Therefore in order to
generalize \cite{BMR} to other interesting situations, one has to
develop a geometric quantization machinery as in \cite{bk1},
\cite{bk3}. One can in fact hope to generalize \cite{BB} as well --
the following has been conjectured in \cite{ka1}.

\begin{conj}
In the assumptions of Theorem~\ref{qua}, the global sections
functor
$$
\Shv(\X,\calo_h) \to \Shv (S[[h]],\pi_*\calo_h)
$$
from sheaves of finitely generated $\calo_h$-modules on $\X$ to
sheaves of finitely generated $\pi_*\calo_h$-modules on $S[[h]] = S
\whtimes \Spec K[[h]]$ is an equivalence of abelian categories over
a dense open subset $U \subset S[[h]]$.
\end{conj}

However, from the present perspective, another relation to Geometric
Representation Theory seems more promising.

\medskip

Let us summarize once more the main steps in the construction of a
tilting generator on a symplectic resolution $X$.
\begin{enumerate}
\item We reduce $X$ to a smooth symplectic variety $X_\m$ over a
  perfect field $k$ of positive characteristic $\cchar k =p > 0$.
\item Using quantization theory, we deform the Frobenius map
$$
\sigma:\calo_{X_\m^\tw} \to \Fr_*\calo_X
$$
to a central algebra map
$$
s:\calo_{X_\m^\tw} \to \calo_h.
$$
\item Using rigidity of tilting vector bundles, we lift the map $s$
  to a central algebra map
\begin{equation}\label{q.fr}
\calo_X \to \calr,
\end{equation}
where $\calr = \eend(\E)$ is a matrix algebra sheaf on $X$.
\end{enumerate}
In principle, a similar procedure can be applied to a smooth variety
$X$ which is only Poisson, not symplectic. The problem is, steps
\thetag{ii} and \thetag{iii} require some rigidity, and one cannot
expect them to work nearly as well for arbitrary Poisson
varieties. In particular, our approach to quantization is
essentially that of Fedosov \cite{fedo}, and it is based on the
fact that locally, all symplectic manifolds and all quantizations
are the same -- this of course breaks down completely in the general
Poisson case.

\medskip

However, there is one more situation where quantization works really
well, namely, the case of a semisimple Lie group with a Poisson-Lie
structure (see \cite{dr}). In this case, the necessary rigidity is
provided by the fact that $G$ is a group -- a quantization becomes
essentially unique if one requires it to be compatible with the
group structure (see \cite{dr} and also \cite{EK}). Motivated by
this, we expect, roughly, the following picture in the Lie group
case.
\begin{quote}
{\em Assume given a semisimple Lie group $G$ over a field $K$ of
characteristic $0$. Then the standard Poisson-Lie structure on $G$
canonically extends to a model $G_O$ of the group $G$ over a
subalgebra $O \subset K$ of finite type over $\Z$, so that for any
maximal ideal $\m \subset O$, we obtain a Poisson-Lie group $G_\m$
over a finite field $k = O/\m$. The Poisson-Lie group $\G_\m$ admits
a unique Frobenius-constant quantization compatible with the group
structure. Moreover, the quantized structure sheaf $\calo_h$ on
$G_\m^\tw$ lifts uniquely to a algebra sheaf $\calo_q$ of finite
rank on $G$ which is, again, compatible with the group structure.}
\end{quote}
This is much too imprecise to be stated even as a Conjecture. In
particular, one has to clarify the exact meaning of ``compatibility
with the group structure'' -- we expect that it should not be
difficult to do this, but at present, this has not been done. In
addition, one cannot expect $\calo_h(h^{-1})$ to be a matrix
algebra, so that step \thetag{iii} -- lifting to characteristic $0$
-- will not be automatic, and probably requires the same methods as
step \thetag{ii}.

In spite of all this, we can guess what the final result will be --
that is, what is the algebra sheaf $\calo_q$. Namely, recall that
G. Lusztig -- see, e.g., \cite{Lu} -- has found a particular form
$U_q$ of the quantized enveloping algebra $U_h$ associated in
\cite{dr} to a semisimple Lie group $G$ (we note that this is {\em
different} from the quantized enveloping algebra used in
\cite{Kobi}). The algebra $U_q$ is defined over a much smaller
subalgebra $K[q,q^{-1}] \subset K[[h]]$ in the algebra of formal
power series in $h = \log q$. Therefore one can actually assign some
value to the parameter $q$. It is known that the resulting algebra
is especially interesting when $q$ is a root of $1$. In this case,
G. Lusztig constructs in addition the so-called {\em quantum
Frobenius map} -- an algebra map $U_q \to U$ from $U_q$ to the usual
universal enveloping algebra $U$ associated to the group $G$. The
dual picture has been also studied, for instance in \cite{CP}. There
instead of quantized enveloping algebra $U_q$, one considers a
quantum version $\calo_{G,q}$ of the algebra $\calo_G$ of algebraic
functions on $G$; if $q$ is a root of $1$, one obtains a quantum
Frobenius map
$$
\calo_G \to \calo_{G,q},
$$
so that $\calo_{G,q}$ becomes a sheaf of associative algebras on the
group $G$.

\medskip

This is what we expect our sheaf $\calo_q$ to be, for $q =
\exp\frac{2\pi\sqrt{-1}}{p}$. The map \eqref{q.fr} should be the
quantum Frobenius map.

\medskip

Unlike \cite{Lu} and consequently \cite{CP}, where $U_q$ and
$\calo_{G,q}$ are constructed by explicit formulas, it should be
possible to obtain $\calo_q$ by pure deformation theory, as an
essentially unique solution to a deformation problem. We do not know
whether it has any real significance for the theory of quantum
groups, a subject very well studied already; still, a conceptual
explanation of the formulas in \cite{Lu} may be worth trying for.

Conversely, the algebras $R = \eend(\E)$ constructed in
Theorem~\ref{equi.thm} should be related to quantum group theory, at
least in the cases like Example~\ref{homo.exa} when the symplectic
resolution $X$ is related to a semisimple group $G$. When $X$ is
{\em not} directly related to any group, one could still hope to
find in the algebras $R$ some of the rich additional structures
known for quantum groups, such as e.g. the so-called {\em crystal
bases}. From this point of view, the most promising case is perhaps
Example~\ref{qui.exa}, the quiver variety case. Since a quiver
variety $X$ is given by a very explicit set of combinatorial data,
computing the algebra $R$ explicitly is not perhaps quite out of
reach.

\bigskip

\noindent
{\sc Steklov Math Institute\\
Moscow, USSR}

\bigskip

\noindent
{\em E-mail address\/}: {\tt kaledin@mccme.ru}

\end{document}